\documentclass[11pt]{amsart}
\usepackage{mathrsfs}
\usepackage{amssymb}

\usepackage{titletoc}
\usepackage{amscd}
\usepackage{amsmath, amssymb}
\usepackage{amsthm}
\usepackage{amsfonts}
\usepackage[colorlinks,linkcolor=blue,citecolor=blue, pdfstartview=FitH]{hyperref}

\usepackage{stmaryrd}
\usepackage{tikz-cd} 
\usepackage{stackrel} 
\usepackage{mathtools}

\numberwithin{equation}{section}

\theoremstyle{plain}
\newtheorem{thm}{Theorem}[section]

\newtheorem{facts}[thm]{Facts}
\newtheorem{quest}[thm]{Question}
\newtheorem{lem}[thm]{Lemma}

\newtheorem{prop}[thm]{Proposition}

\theoremstyle{definition}

\newtheorem{rmk}[thm]{Remark}

\newtheorem{definition}[thm]{Definition}

\newtheorem{exm}[thm]{Example}

\newtheorem{defn-thm}[thm]{Definition-Theorem}
\newtheorem{defn-pro}[thm]{Definition-Proposition}
\newtheorem{conjecture}[thm]{Conjecture}

\usepackage{fancyhdr}

\newcommand{\A}{{\mathbb A}}

\newcommand{\D}{{\mathbb D}}

\newcommand{\sE}{{\mathcal E}}
\newcommand{\scE}{{\mathscr E}}
\newcommand{\F}{{\mathbb F}}

\newcommand{\sF}{{\mathcal F}}
\newcommand{\scF}{{\mathscr F}}

\newcommand{\G}{{\mathbb G}}
\newcommand{\sG}{{\mathcal G}}
\newcommand{\scG}{{\mathscr G}}

\newcommand{\sH}{{\mathcal H}}

\newcommand{\sI}{{\mathcal I}}

\renewcommand{\L}{{\mathbb L}}

\newcommand{\scL}{{\mathscr L}}

\newcommand{\sO}{{\mathcal O}}
\renewcommand{\P}{{\mathbb P}}

\newcommand{\sP}{{\mathcal P}}
\newcommand{\scP}{{\mathscr P}}
\newcommand{\Q}{{\mathbb Q}}

\newcommand{\R}{{\mathbb R}}

\newcommand{\sT}{{\mathcal T}}
\newcommand{\scT}{{\mathscr T}}

\newcommand{\sW}{{\mathcal W}}
\newcommand{\X}{{\mathbb X}}

\newcommand{\Z}{{\mathbb Z}}

\renewcommand{\k}{\Bbbk} 
\newcommand{\id}{\mathrm{id}}

\newcommand{\Hom}{\mathrm{Hom}}
\newcommand{\End}{\mathrm{End}}

\newcommand{\Ext}{\mathrm{Ext}}
\newcommand{\op}{{\mathrm{op}}}
\newcommand{\rk}{{\mathrm{rk}}}

\newcommand{\Rep}{\mathrm{Rep}}

\newcommand{\Perv}{\mathrm{Perv}} 

\newcommand{\Ind}{\mathrm{Ind}} 
\newcommand{\Res}{\mathrm{Res}} 
\newcommand{\Sch}{\mathrm{Sch}}
\newcommand{\Sets}{\mathrm{Sets}}
\newcommand{\AffSch}{\mathrm{AffSch}}
\newcommand{\Alg}{\mathrm{Alg}}
\newcommand{\Dbc}{D^{\mathrm{b}}_{\mathrm{c}}}
\newcommand{\Db}{D^{\mathrm{b}}}

\newcommand{\Spec}{\mathrm{Spec}}

\newcommand{\IC}{\mathrm{IC}} 
\newcommand{\af}{\mathrm{aff}}

\newcommand{\cLim}[1]{\mathop{\mathrm{colim}}\limits_{#1}}
\newcommand{\br}[1]{[ \hspace{-1pt} [ #1 ] \hspace{-1pt} ]}
\newcommand{\pa}[1]{( \hspace{-1pt} ( #1 ) \hspace{-1pt} )}

\newcommand{\Gr}{\mathrm{Gr}} 
\newcommand{\Conv}{\mathrm{Conv}} 
\newcommand{\GGr}{\mathbb{G}\mathrm{r}} 
\newcommand{\Fl}{\mathrm{Fl}} 

\newcommand{\Gm}{{\mathbb{G}_{\mathrm{m}}}}
\newcommand{\Loop}{\mathrm{L}}
\newcommand{\simto}{\xrightarrow{\sim}}
\mathchardef\mhyphen="2D
\newcommand{\rS}{\mathrm{S}}
\newcommand{\rT}{\mathrm{T}}
\newcommand{\pH}{{}^{\mathrm{p}} \hspace{-2pt} \mathscr{H}}
\newcommand{\coH}{\mathsf{H}}
\newcommand{\coHc}{\mathsf{H}_{\mathrm{c}}}

\newcommand{\tboxtimes}{\mathbin{\widetilde{\mathord{\boxtimes}}}}
\newcommand{\coweyl}{\mathsf{N}}
\newcommand{\weyl}{\mathsf{M}}
\newcommand{\simp}{\mathsf{L}}
\newcommand{\til}{\mathsf{T}}
\newcommand{\un}{\mathrm{u}}

\makeatletter
\def\lotimes{\@ifnextchar_{\@lotimessub}{\@lotimesnosub}}
\def\@lotimessub_#1{\mathchoice{\mathbin{\mathop{\otimes}^L}_{#1}}%
  {\otimes^L_{#1}}{\otimes^L_{#1}}{\otimes^L_{#1}}}
\def\@lotimesnosub{\mathbin{\mathop{\otimes}^L}}

\def\lboxtimes{\@ifnextchar_{\@lboxtimessub}{\@lboxtimesnosub}}
\def\@lboxtimessub_#1{\mathchoice{\mathbin{\mathop{\boxtimes}^L}_{#1}}%
  {\boxtimes^L_{#1}}{\boxtimes^L_{#1}}{\boxtimes^L_{#1}}}
\def\@lboxtimesnosub{\mathbin{\mathop{\boxtimes}^L}}
\makeatother

\begin{document}
\title{Some applications of the geometric Satake equivalence to modular representation theory} 

\makeatletter
\let\MakeUppercase\relax
\makeatother


\author{Simon Riche}
\address{Universit\'e Clermont Auvergne, CNRS, LMBP, F-63000 Clermont-Ferrand, France.}
\email{simon.riche@uca.fr}



\begin{abstract}
These notes present an application of the geometric Satake equivalence to the description of characters of indecomposable tilting modules for reductive algebraic groups over fields of positive characteristic, obtained in joint work with G. Williamson.
\end{abstract}

\maketitle 


\section*{Introduction}

\subsection{Geometric Satake equivalence and Representation Theory}

The geometric Satake equivalence is an equivalence of categories relating the category of perverse sheaves on the affine Grassmannian of a connected reductive algebraic group and the category of representations of the Langlands dual group. This equivalence was initially considered (following ideas of Lusztig and Drinfeld) in the case of sheaves with coefficients in a field of characteristic $0$, and found in this setting celebrated applications in various directions, including the geometric Langlands program and representation theory (in particular via the notion of Mirkovi\'c--Vilonen cycles). However, the construction of this equivalence by Mirkovi\'c--Vilonen~\cite{MV07} applies to more general coefficients, in particular fields of positive characteristic. In this setting, applications of this result were developed only recently. The goal of these notes is to explain such an application in representation theory, found in recent joint work with Geordie Williamson~\cite{RW22}. More explicitly, using Smith--Treumann theory we were able to extract from the geometric Satake equivalence a character formula for indecomposable tilting modules in all blocks of representations of a reductive algebraic group over a field of arbitrary characteristic $\ell$, in terms of the $\ell$-canonical basis~\cite{gjw}. This solved a conjecture of ours formulated in~\cite{RW18}.

\subsection{Contents}


Sections~\ref{sec:lecture-1} and~\ref{sec:lecture-2} are devoted to a brief presentation of the geometric Satake equivalence.
There already exist a number of reviews of this subject of various lengths and degrees of technicality, among which~\cite{Zhu16,BR18,workshop} and~\cite[\S 1]{central}. We do not feel it is worth writing another such account; instead we have tried to be as direct as possible, getting as fast as possible to the main statements, and discussing in detail only some important technical constructions which are sometimes omitted in the standard references on this subject.

Sections~\ref{sec:lecture-3} and~\ref{sec:lecture-4} are concerned with the application of this construction to modular representation theory of reductive algebraic groups. Here we discuss the study of perverse sheaves corresponding to tilting modules under the geometric Satake equivalence (which involves the ``Iwahori--Whittaker model'' for the Satake category studied in joint work with Bezrukavnikov, Gaitsgory, Mirkovi{\'c} and Rider~~\cite{BGMRR19}), and the use of Smith--Treumann theory to ``extract'' combinatorial information using this realization, obtained in~\cite{RW22}.

\subsection{Some notation}

We will use the following (standard) notations for categories:
\begin{itemize}
\item
$\Sets$: category of sets;
\item
$\Sch_S$: category of $S$-schemes (for $S$ a base scheme);
\item
$\Alg_R$: category of commutative $R$-algebras (for $R$ a fixed commutative ring);
\item
$\AffSch_R$: category of affine schemes over $\Spec(R)$ (for $R$ a fixed commutative ring).
\end{itemize}
If $R$ is a commutative ring, we will also write $\Sch_R$ for $\Sch_{\Spec(R)}$.

\subsection{Acknowledgements}

This text grew up of the notes prepared for the mini-course \emph{Applications of the geometric Satake correspondence} given at the conference \emph{Th\'eorie des representations \`a Lyon} in June 2023, which were typed by Quan Situ.
We thank him for offering to type these notes, which offered an indispensable input for this project, and Florence Fauquant-Millet and Philippe Gille for organizing this event and offerring me the opportunity to give these lectures.

This project has received
funding from the European Research Council (ERC) under the European Union's Horizon 2020
research and innovation programme (grant agreements No~101002592).

\section{Technical preliminaries} 
\label{sec:lecture-1}

In this section we recall some general definitions and results regarding ind-schemes (in~\S\ref{ss:ind-schemes}), affine Grassmannians (in~\S\ref{ss:Gr}), and attractors for actions of the multiplicative group $\Gm$ on schemes (in~\S\ref{ss:attractors}). 
Our main references are~\cite{Richarz20} for the first two topics, and~\cite{Richarz16} for the third one. (All of these results have earlier variants, but often with stronger assumptions. See~\cite{Richarz20, Richarz16, Zhu16} for references.)

\subsection{Ind-schemes} 
\label{ss:ind-schemes}

Let $\k$ be a field. 
Consider the category $\Sch_\k$ of $\k$-schemes, and the subcategory $\AffSch_\k$ of affine schemes, which is equivalent to the opposite category of the category $\Alg_\k$. 
Recall that to any $\k$-scheme $X$ one can attach its functor(s) of points
\[
h_X=\Hom_{\Sch_{\Spec(\k)}}(-,X) \colon \Sch_\k^\op \rightarrow \Sets, \quad
h'_X=\Hom_{\Sch_{\Spec(\k)}}(-,X) \colon \AffSch_\k^\op \rightarrow \Sets.
\]
By the Yoneda lemma the assignment $X \mapsto h_X$ is fully-faithful, and the fact that any scheme admits an affine open cover implies that the assignment $X \mapsto h'_X$ is also fully-faithful.
In other words, a scheme is determined by each version of its functor of points. For simplicity, we will say that a functor ``is a scheme'' if it is representable by a scheme, i.e.~of the form $h_X$ (or $h'_X$) for some $\k$-scheme $X$, and in this case we will write $X$ for $h_X$ (or $h_X'$).

\begin{rmk}
\label{rmk:Sch-AffSch}
Depending on the contexts, one might want to work with the functors $h_X$ (defined on $\Sch_\k^\op$) or $h_X'$ (defined on $\AffSch_\k^\op$). The latter one seems to be the best choice in the setting we consider here, as illustrated in~\cite{Richarz20} (and below). Note that, as explained in~\cite[Remark~1.6]{Richarz20}, any functor $\AffSch_\k^\op \rightarrow \Sets$ which is a Zariski sheaf extends in a unique way to a functor $\Sch_\k^\op \rightarrow \Sets$ which is a Zariski sheaf.
\end{rmk}

\begin{definition}
An ind-scheme over $\k$ is a functor 
\[
X \colon \AffSch_\k^\op \rightarrow \Sets
\]
which admits a presentation $X= \cLim{i} X_i$ for some filtered poset $(I, \leq)$ and some family of schemes $(X_i : i \in I)$, where each transition morphism $h_{X_i} \rightarrow h_{X_j}$ ($i \leq j$) is (induced by) a closed immersion of schemes.
\end{definition}

Note that colimits of functors are computed pointwise, so that the fact that $X= \cLim{i} X_i$ means that for any $S \in \AffSch_\k$ we have
\[
X(S) = \cLim{i} \Hom(S,X_i).
\]

If $X$, $Y$ are ind-schemes, a morphism of ind-schemes $X \to Y$ is simply a morphism of functors from $X$ to $Y$.
Ind-schemes in this sense are sometimes called \emph{strict} ind-schemes. Since all the ind-schemes we want to consider will be of this form, we will omit the adjective ``strict.'' 
Given an ind-scheme $X$, the datum of a poset $I$, a family of scheme $(X_i : i \in I)$ and transition morphisms $X_i \to X_j$ such that $X = \cLim{i} X_i$ is called a \emph{presentation} of $X$.

\begin{exm}
\phantomsection
\label{ex:ind-schemes}
\begin{enumerate}
	\item Each scheme naturally defines an ind-scheme (with the set ``$I$'' having only one element). 
	\item 
	\label{it:ind-sch-affine-space}
	For any set $I$, the functor
	\[
	\A^I_\k \colon \left\{ 
	\begin{array}{ccc}
	\AffSch_\k^\op & \rightarrow & \Sets \\
	Y & \mapsto & \bigoplus\limits_{i\in I} \Gamma(Y,\sO_Y)
	\end{array}
	\right.
	\]
	is an ind-scheme. 
	Indeed, we have $\A^I_\k=\cLim{J\subset I} \A_\k^J$, where $J$ runs over all finite subsets of $I$. 
\end{enumerate}
\end{exm}

\begin{rmk}
Regarding Example~\ref{ex:ind-schemes}\eqref{it:ind-sch-affine-space}, in case $I$ is finite the ind-scheme $\A_\k^I$ is of course a scheme (an affine space of dimension $\# I$). But if $I$ is infinite $\A_\k^I$ is \emph{not} a scheme. In fact, assume for a contradiction that it were, and denote this scheme by $X$. Then there exists a ring $A$ and an open immersion $\Spec(A) \to X$. By definition there exists a finite subset $J \subset I$ such that $j$ factors through a morphism $\Spec(A) \to \A_\k^J$. Then, for any finite subset $K \subset I$ containing $J$ we have
\[
\Spec(A) = \Spec(A) \times_{X} \A_\k^K,
\]
hence we have an open immersion $\Spec(A) \to \A_\k^K$. This implies that $\dim(A) = \# K$ for any such $K$ (e.g.~by~\cite[Theorem~5.22(3)]{goertz-wedhorn}), which is absurd.


In particular, one should not confuse the ind-scheme $\A_\k^I$ with the affine space given by $\Spec(\k[X_i : i \in I])$.
\end{rmk}

\begin{rmk} 
\phantomsection
\label{rmk:ind-schemes}
\begin{enumerate}
	\item 
	\label{it:ind-sch-fpqc}
	An ind-scheme is automatically a sheaf for the fpqc topology; see~\cite[Lemma~1.4]{Richarz20}.
	\item Let $X$ be an ind-scheme, and consider a presentation $X=\cLim{i} X_i$.
For any scheme $Y$ we have a natural morphism 
\[
\cLim{i} \Hom(Y, X_i) \rightarrow \Hom(Y,X).
\]
This morphism is always injective, but not always a bijection; it is so however if $Y$ is quasi-compact. (This easily follows from the similar claim when $Y$ is affine, which is a tautology.) As a consequence, if $X$, $Y$ are ind-schemes with respective presentations $X=\cLim{i} X_i$, $Y=\cLim{j} Y_i$, and if each $X_i$ is quasi-compact, we have
\begin{equation}
\label{eqn:Hom-ind-sch}
\Hom(X,Y) = \lim_i \cLim{j} \Hom(X_i, Y_j).
\end{equation}
\end{enumerate}
\end{rmk}

By definition, ind-schemes are functors from $\AffSch_\k^\op$ to $\Sets$. One can form in the obvious way fiber products of functors. It turns out that given ind-schemes $X$, $Y$, $Z$ and morphisms $X \to Z$, $Y \to Z$, the fiber product $X \times_Z Y$ is an ind-scheme; for details see~\cite[Lemma~1.10(1)]{Richarz20}.

\begin{definition}
Let $X,Y$ be functors from $\AffSch_\k^\op$ to $\Sets$
(e.g.~ind-schemes). 
A morphism of functors $f:X\rightarrow Y$ is said to be \emph{representable by a closed}, resp.~\emph{open}, resp.~\emph{locally-closed immersion}, if for any scheme $Z$ and any morphism $Z\rightarrow Y$ the product $X\times_Y Z$ is a scheme and the induced morphism $X\times_Y Z\rightarrow Z$ is a closed, resp.~open, resp.~locally-closed, immersion of schemes. 
\end{definition}

\begin{definition}
Let $X$ be an ind-scheme. 
We will say that $X$ is \emph{separated}, resp.~\emph{ind-affine}, resp.~\emph{of ind-finite type}, resp.~\emph{ind-proper}, if it admits a presentation $X=\cLim{i} X_i$ where each $X_i$ is separated, resp.~affine, resp.~of finite type, resp.~proper. 
\end{definition} 

For an explanation why we do not add ``ind'' to the adjective ``separated,'' see~\cite[Exercise~1.31]{Richarz20}.

\subsection{Affine Grassmannians}
\label{ss:Gr}

\subsubsection{Definition}

We start by recalling the definition of a torsor for a group scheme. For a base scheme $S$, we denote by $\mathrm{Sch}_S$ the category of $S$-schemes.

\begin{definition}
\label{def:torsors}
Let $S$ be a scheme, and consider a group scheme $G\rightarrow S$. 
A \emph{$G$-torsor} is a fppf sheaf $\sP$ on the category $\mathrm{Sch}_S$ endowed with a right action of $G$ such that:
\begin{enumerate}
	\item 
	\label{it:covering-torsors}
	for any $X\in \mathrm{Sch}_S$, there exists a fppf covering $(X_i\rightarrow X)_{i\in I}$ such that $\sP(X_i)\neq \varnothing$ for all $i$; 
	\item the map $\sP \times_S G \rightarrow \sP \times_S \sP$ defined by $(x,g)\mapsto (x,x \cdot g)$ is an isomorphism of sheaves. 
\end{enumerate}
\end{definition} 

Any group scheme $G \to S$ has a \emph{trivial torsor}, namely $G$ itself (for the action given by mutliplication on the right).
We note the following properties (for details, see~\cite[\S 2.1.2]{central}):
\begin{itemize}
\item
given a morphism of schemes $T\rightarrow S$, for any $G$-torsor $\sP$ we have an associated $G\times_S T$-torsor $\sP \times_S T$;
\item 
if $G\rightarrow S$ is affine, then any $G$-torsor is represented by a scheme over $S$;
\item
if moreover $G\rightarrow S$ is smooth, then any $G$-torsor is \'{e}tale locally trivial. (In other words, the covering in Item~\eqref{it:covering-torsors} of Definition~\ref{def:torsors} can be chosen to be an \emph{\'etale} covering.)
\end{itemize}

From now on we fix a field $\k$.
For $R\in \Alg_\k$, we set 
\[
\D_R=\Spec (R \br{z}), \qquad \D_R^\times=\Spec (R \pa{z} ).
\]
Given an algebra morphism $R \to S$ we have induced algebra morphisms $R\br{z} \to S\br{z}$ and $R\pa{z} \to S\pa{z}$, hence morphisms of schemes
\[
\D_S \to \D_R, \qquad \D_S^\times \to \D_R^\times.
\]

We now consider a
flat affine group scheme $G\rightarrow \D_\k$ of finite type. 

\begin{definition}
The \emph{affine Grassmannian for $G$} is the functor 
\[
\GGr_G \colon \Alg_\k \rightarrow \Sets
\]
which sends $R\in \Alg_\k$ to the set of isomorphism classes of pairs $(\sE, \alpha)$ where $\sE\rightarrow \D_R$ is a $(G \times_{\D_\k} \D_R)$-torsor and $\alpha\colon \D_R^\times \to \sE$ is a section of $\sE$ over $\D_R^\times$ (i.e.~a morphism whose composition with the projection $\sE\rightarrow \D_R$ is the natural embedding $\D_R^\times \hookrightarrow \D_R$).
\end{definition} 

\begin{rmk} 
Given a bundle $\sE$, the datum of the section $\alpha$ is equivalent to the datum of an isomorphism (called a \emph{trivialization}) 
\[
\sE_{|\D_R^\times}\simto (G \times_{\D_\k} \D_R) \times_{\D_R} \D_R^\times = G\times_{\D_\k} \D_R^\times.
\]
\end{rmk}

An important special case of this construction is when $G=H \times_{\Spec(\k)} \D_\k$ for some group scheme $H$ over $\k$. We introduce a special notation for this case, setting
\[
\Gr_H=\GGr_{H \times_{\Spec (\k)} \D_\k}.
\] 

\subsubsection{Representablity}

We continue with our flat affine group scheme $G\rightarrow \D_\k$ of finite type.
The following result is~\cite[Theorem~3.4]{Richarz20}.

\begin{thm}
\phantomsection
\label{thm:GrG-indscheme}
\begin{enumerate}
\item The affine Grassmannian $\GGr_G$ is represented by a separated ind-scheme of ind-finite type over $\k$. 
\item If $G$ is reductive, then $\GGr_G$ is ind-projective. 
\end{enumerate}
\end{thm}

The idea of the proof of this theorem goes as follows. First, one checks the statement directly in case $G=\mathrm{GL}_{n,\k \br{z}}$, by describing the affine Grassmannian as a moduli space of lattices. Then one 
shows that if $H\rightarrow G$ is a closed immersion of flat affine group schemes of finite type such that the fppf quotient $G/H$ is representable by a quasi-affine (resp.~affine) scheme, then the induced morphism $\GGr_H\rightarrow \GGr_{G}$ is representable by a quasi-compact immersion (resp.~a closed immersion). 
Finally, one shows that for any $G$ as in the theorem there exists $n \geq 1$ and a closed immersion of group schemes $G \rightarrow \mathrm{GL}_{n,\k \br{z}}$ such that the fppf quotient $\mathrm{GL}_{n,\k \br{z}} / G$ is representable by a quasi-affine scheme, which is automatically affine in case $G$ is reductive. 

The idea of this proof goes back at least to~\cite[\S 4.5]{beilinson-drinfeld}. This idea has been repeated in various degrees of generality since; see e.g.~\cite[\S 3]{heinloth} or~\cite[\S 1.2]{Zhu16}.

\subsubsection{Beauville--Laszlo gluing} 

Let $X\rightarrow \Spec (\k)$ be a scheme, and let $x\in X$ be a $\k$-point of $X$ such that the local ring $\sO_{X,x}$ of $X$ at $x$ is regular of dimension $1$. (In practice $X$ will be chosen to be a smooth curve over $\k$; in fact the case when $X=\A^1_\k$ suffices.) In this setting, as explained in~\cite[\S 3.2]{Richarz20}, the completion of the local ring $\sO_{X,x}$ is isomorphic (as a $\k$-algebra) to $ \k \br{z}$. We fix such an isomorphism; this provides a morphism of schemes $\D_\k \to X$.

Let $\sG\rightarrow X$ be a flat affine group scheme of finite presentation. 
Consider the functor 
\[
\GGr_{\sG,x}\colon \Alg_\k\rightarrow \Sets
\]
sending $R\in \Alg_\k$ to the set of isomorphism classes of pairs $(\sE,\alpha)$ where $\sE$ is a $\sG \times_{\Spec (\k)} \Spec (R)$-torsor 
and
\[
\alpha \colon (X \smallsetminus \{x\}) \times_{\Spec(\k)} \Spec (R) \rightarrow \sE
\]
is a section of $\sE$ over $(X \smallsetminus \{x\}) \times_{\Spec(\k)} \Spec (R)$.

The following result is often called the Beauville--Laszlo gluing theorem (because the first version of such a statement is due to Beauville--Laszlo, see~\cite{beauville-laszlo}). See~\cite[Theorem~3.15]{Richarz20} or~\cite[Theorem~1.4.3]{Zhu16}.

\begin{thm} 
\label{thm:beauville-laszlo}
Set $G=\sG\times_{X} \D_\k$. Then
there exists a canonical isomorphism 
\[
\GGr_{\sG,x}\simto \GGr_G.
\] 
\end{thm}

This theorem is important because it allows to pass from the ``local'' problem described by $\GGr_G$ (in terms of torsors over the infinitesimal disc $\D_\k$) to a ``global'' problem described in terms of torsors over the scheme $X$.

\subsubsection{Loop spaces} 

For a functor $X\colon \Alg_{\k\pa{z}}\rightarrow \Sets$, we define the associated \emph{loop functor} 
\[
\L X\colon \Alg_{\k} \rightarrow \Sets, \quad R \mapsto X\big(R \pa{z}\big).
\]
For a functor $X\colon \Alg_{\k\br{z}}\rightarrow \Sets$, we define the associated \emph{arc functor} 
\[
\L^+ X\colon \Alg_{\k} \rightarrow \Sets, \quad R \mapsto X\big(R\br{z}\big).
\] 
For $i \geq 0$ and a functor $X\colon \mathrm{Aff}_{\k[z]/(z^{i+1})}\rightarrow \Sets$, we define the associated \emph{jet functor} 
\[
\L^+_i X\colon \Alg_{\k} \rightarrow \Sets, \quad R \mapsto X\big(R[z]/(z^{i+1})\big).
\]

We will mainly consider these functors in case the functor is constructed from a scheme over $\D_\k$. In this case we will simplify the notation as follows: for a scheme $X$ over $\D_\k$ and $i \in \Z_{\geq 0}$ we write 
\[
\L X:=\L\big(X\times_{\D_\k} \D_\k^\times \big), \quad
\L^+_i X=\L^+_i\big(X\times_{\D_\k} \Spec (\k[z]/(z^{i+1}))\big).
\]
In this case, if $i \geq j$ we have a natural ``reduction'' morphism $\L^+_j X \to \L^+_i X$.
An even more specific case is when the scheme $X$ is obtained by base change from a $\k$-scheme $Y$. In this case we set
\[
\Loop Y=\L\big(Y \times_{\Spec(\k)} \D_\k \big), \quad 
\Loop^+ Y=\L\big(Y\times_{\Spec (\k)} \D_\k \big), \quad
\Loop^+_i Y = \Loop^+_i (Y \times_{\Spec(\k)} \D_\k). 
\]
The functor $\Loop Y$ should be thought of as ``infinitesimal loops in $Y$'', i.e.~maps from the punctured infinitesimal disc $\D_\k^\times$ to $Y$. A similar mental picture for $\Loop^+ Y$ can be obtained by replacing infinitesimal loops by infinitesimal arcs, i.e.~maps from the infinitesimal disc $\D_\k$ to $Y$.

The following statement is standard; see~\cite[Lemma~3.17]{Richarz20} and~\cite[\S 1.1]{cesnavicius}.

\begin{prop}
\phantomsection
\label{prop:loop-functors}
\begin{enumerate}
\item Let $X$ be a scheme over $\D_\k$. 
Then each $\L^+_i X$ is a scheme over $\k$, which is smooth, resp.~of finite type, if $X$ is, each reduction morphism $\L^+_j X\rightarrow \L^+_i X$ is an affine morphism, and
\[
\L^+ X \simeq \varprojlim_{i} \L^+_i X
\]
is a scheme. 
If $X$ is affine, then $\L^+ X$ is also affine. 
\item If $X\rightarrow \D_\k^\times$ is an affine scheme, then $\L X$ is an ind-affine ind-scheme. 
\end{enumerate}
\end{prop}

\subsubsection{Uniformization}

Let $G\rightarrow \D_\k$ be a smooth affine group scheme of finite type. Then $\L G$ and $\L^+ G$ naturally factor through group-valued functors.
In fact, by Proposition~\ref{prop:loop-functors}, $\L G$ is an ind-affine group ind-scheme over $\k$, and $\L^+ G$ is an affine group scheme over $\k$. We also have a canonical inclusion $\L^+ G\subset \L G$. We can therefore consider the presheaf quotient $\L G / \L^+ G$, i.e.~the functor on $\Alg_\k$ sending $R$ to the quotient group $(\L G)(R) / (\L^+ G)(R)$. We will denote by $(\L G/ \L^+ G )_{\mathrm{\acute{e}t}}$ the \'etale sheafification of this presheaf.

For the following statement we refer to~\cite[Proposition~3.18]{Richarz20}.

\begin{thm}
\label{thm:Gr-quotient}
There exists a canonical isomorphism of \'{e}tale sheaves 
\[
\big(\L G/ \L^+ G\big)_{\mathrm{\acute{e}t}} \simto \GGr_G
\]
induced by $g\mapsto (\sE_0, g^{-1})$, where $\sE_0$ is the trivial torsor. 
\end{thm}

\begin{rmk} 
\begin{enumerate}
	\item This theorem is closely related to the property that if $\sE\rightarrow \D_R$ is a $(G \times_{\D_\k} \D_R)$-torsor then there exists an \'{e}tale cover $R\rightarrow R'$ such that $\sE\times_{\D_R} \D_{R'}$ is isomorphic to the trivial $(G \times_{\D_\k} \D_{R'})$-torsor. 
	\item In view of Remark~\ref{rmk:ind-schemes}\eqref{it:ind-sch-fpqc}, this theorem shows in particular that the \'etale sheaf $(\L G/ \L^+ G)_{\mathrm{\acute{e}t}}$ is in fact an fpqc sheaf. One can say more in the (important!) special case when $G=H \times_{\Spec(\k)} \D_\k$ for some connected reductive group $H$ over $\k$: in this case, the Zariski sheafification $(\Loop H/ \Loop^+ H)_{\mathrm{Zar}}$ coincides with $(\Loop H/ \Loop^+ H)_{\mathrm{\acute{e}t}}$ (hence is a fpqc sheaf), see~\cite[Theorem~2.5]{cesnavicius}, and if $H$ is split over $\k$ (e.g.~if $\k$ is algebraically closed) then the presheaf quotient $\Loop H/ \Loop^+ H$ coincides with $(\Loop H/ \Loop^+ H)_{\mathrm{\acute{e}t}}$, see~\cite[Theorem~3.4 and Example~3.2]{cesnavicius}.
\end{enumerate} 
\end{rmk} 

\subsection{Attractors and Braden's theorem} 
\label{ss:attractors}


\subsubsection{Definition} 

Let $X\rightarrow \Spec (\k)$ be a scheme endowed with an action of the multiplicative group $\G_{\mathrm{m},\k}$. For a $\k$-scheme $T$ we set $\G_{\mathrm{m},T} = \G_{\mathrm{m},\k} \times_{\Spec(\k)} T$. We also denote by $(\A^1_T)^+$, resp.~$(\A^1_T)^-$, the scheme given by the affine line over $T$, endowed with the natural action of $\G_{\mathrm{m},T}$, resp.~the opposite of the natural action.

We define the attractor, repeller, and fixed points of the action as follows.

\begin{definition} 
\phantomsection
\label{def:attractors}
\begin{enumerate}
\item
The functor of fixed points $X^\circ$ is
the functor 
\[
X^\circ \colon \AffSch_\k^\op \rightarrow \Sets
\] 
sending a $\k$-scheme $T$ to the set of $T$-morphisms $T\rightarrow X\times_{\Spec(\k)} T$ such that the diagram
\[
\begin{tikzcd}
\G_{\mathrm{m},T} \arrow[r] \arrow[d] & \G_{\mathrm{m},T}\times_T (X\times_{\Spec(\k)} T) \arrow[d] \\ 
T \arrow[r] & X\times_{\Spec(\k)} T
\end{tikzcd}
\] 
commutes, where the right-hand side is induced by the action morphism
\[
\G_{\mathrm{m},\k} \times_{\Spec(\k)} X \to X.
\]
\item
The attractor $X^+$ is the functor
\[
X^+ \colon \AffSch_\k^\op \rightarrow \Sets
\] 
sending $T$ to the sets of $T$-morphisms $(\A^1_T)^+ \rightarrow X\times_{\Spec(\k)} T$ such that the diagram
\[
\begin{tikzcd}
\G_{\mathrm{m},T}\times_T (\A^1_T)^+ \arrow[r] \arrow[d] & \G_{\mathrm{m},T}\times_T (X\times_{\Spec(\k)} T) \arrow[d] \\ 
(\A^1_T)^+ \arrow[r] & X\times_{\Spec(\k)} T
\end{tikzcd}
\]
commutes, where the vertical arrows are induced by the actions of $\G_{\mathrm{m},T}$.
\item
The repeller $X^-$ is the functor
\[
X^- \colon \AffSch_\k^\op \rightarrow \Sets
\] 
sending $T$ to the sets of $T$-morphisms $(\A^1_T)^- \rightarrow X\times_{\Spec(\k)} T$ such that the diagram
\[
\begin{tikzcd}
\G_{\mathrm{m},T}\times_T (\A^1_T)^- \arrow[r] \arrow[d] & \G_{\mathrm{m},T}\times_T (X\times_{\Spec(\k)} T) \arrow[d] \\ 
(\A^1_T)^- \arrow[r] & X\times_{\Spec(\k)} T
\end{tikzcd}
\]
commutes, where the vertical arrows are induced by the actions of $\G_{\mathrm{m},T}$.
\end{enumerate}
\end{definition}

\begin{rmk}
\phantomsection
\begin{enumerate}
\item
Clearly the definitions of $X^\circ(T)$, $X^+(T)$, $X^-(T)$ make sense for any $\k$-scheme $T$, even if $T$ is not affine. These ``extended'' functors are considered in~\cite{Richarz16}, rather than those in Definition~\ref{def:attractors}. This does not lead to any difference in the theory; in fact $X^\circ$, $X^+$ and $X^-$ as defined above are clearly Zariski sheaves, and the versions considered in~\cite[Definition~1.3]{Richarz16} are the unique extensions of these functors to Zariski sheaves from $\Sch_\k$ to $\Sets$ mentioned in Remark~\ref{rmk:Sch-AffSch}. In particular, if $X^\circ=h'_Y$ for some $\k$-scheme $Y$, then the corresponding functor in~\cite{Richarz16} is $h_Y$, and vice versa. Similar considerations apply to $X^+$ and $X^-$.
\item
The definitions above can be considered also in the more general setting of algebraic spaces, as is done in~\cite{Richarz16}. We will restrict to schemes here, since this is sufficient for the applications we want to consider. One can also work over a base scheme $S$ instead of $\Spec(\k)$ (as is also done in~\cite{Richarz16}); this is important in some contexts related to the geometric Satake equivalence, which will however not be discussed here.
\end{enumerate}
\end{rmk}

It is clear that if $Y$ is the scheme $X$ endowed with the opposite action of $\G_{\mathrm{m},\k}$ we have identifications $Y^\circ=X^\circ$, $Y^+=X^-$, $Y^-=X^+$. In particular, the study of repellers can therefore be reduced to that of attractors.

We have natural morphisms of functors
\begin{equation}
\label{eqn:morphisms-att-fixed-pts}
q^+ \colon X^+ \to X^\circ, \quad q^- \colon X^- \to X^\circ
\end{equation}
induced by evaluation at $0 \in \A^1_\k(\k)$, natural morphisms of functors
\begin{equation}
\label{eqn:morphisms-att}
p^+ \colon X^+ \to X, \quad p^- \colon X^- \to X
\end{equation}
induced by evaluation at $1 \in \A^1_\k(\k)$, and a natural morphism of functors
\begin{equation}
\label{eqn:morphism-fixed-pts}
X^\circ \to X.
\end{equation}
There are also natural actions of the functor-in-groups $\G_{\mathrm{m},\k}$ on $X^\circ$, $X^+$ and $X^-$ such that the morphisms in~\eqref{eqn:morphisms-att-fixed-pts},~\eqref{eqn:morphisms-att} and~\eqref{eqn:morphism-fixed-pts} are equivariant.

\begin{rmk}
\label{rmk:att-monomorphism-separated}
Recall that a morphism of schemes $Y \to Z$ is called a \emph{monomorphism} if it is a monomorphism in the category of schemes, i.e.~for any scheme $T$ the morphism $\Hom(T,Y) \to \Hom(T,Z)$ is injective. If $Y$ and $Z$ are $\k$-schemes, then it is equivalent to require that for any $T \in \Sch_\k$ the morphism $\Hom_{\Sch_\k}(T,Y) \to \Hom_{\Sch_\k}(T,Z)$ is injective; in fact this follows from the characterization of monomorphisms in terms of fiber products (see~\cite[\href{https://stacks.math.columbia.edu/tag/08LR}{Tag 08LR}]{stacks-project}) and the fact that the fiber product of $Y$ and $Z$ in $\Sch_\k$ is the same as in the category of all schemes. Using the fact that any scheme admits an affine open cover, one sees that it is equivalent to require that for any $T \in \AffSch_\k$ the morphism $\Hom_{\Sch_\k}(T,Y) \to \Hom_{\Sch_\k}(T,Z)$ is injective.


If $X$ is a \emph{separated} $\k$-scheme, then the morphisms~\eqref{eqn:morphisms-att} are monomorphisms. In fact, by the characterization explained above it suffices to show that for any affine $\k$-scheme $T$ the morphism $\Hom_{\Sch_\k}(T,X^+) \to \Hom_{\Sch_\k}(T,X)$ is injective. Now elements of $\Hom_{\Sch_\k}(T,X^+)$ are in particular morphisms of $T$-schemes $\A^1_T \to X \times_{\Spec(\k)} T$, and two elements whose images in $\Hom_{\Sch_\k}(T,X)$ are equal coincide on the open subscheme $\G_{\mathrm{m}, T} \subset \A^1_T$. Now the scheme theoretic closure of $\G_{\mathrm{m}, T}$ in $\A^1_T$ is $\A^1_T$ (see~\cite[\href{https://stacks.math.columbia.edu/tag/056C}{Tag 056C}]{stacks-project}), which allows to conclude using~\cite[\href{https://stacks.math.columbia.edu/tag/01RH}{Tag 01RH}]{stacks-project}.
\end{rmk}

\subsubsection{Representability} 

The following technical conditions on the action (introduced in~\cite{Richarz16}) turn out to be extremely useful in the study of the functors $X^\circ$, $X^+$ and $X^-$.

\begin{definition} 
The $\G_{\mathrm{m},\k}$-action on $X$ is said to be \emph{\'{e}tale} (resp.~\emph{Zariski}) \emph{locally linearizable} if there exists a $\G_{\mathrm{m},\k}$-equivariant cover $(U_i\rightarrow X)_{i\in I}$ where each $U_i$ is affine and the maps $U_i\rightarrow X$ are \'{e}tale (resp.~open embeddings). 
\end{definition}

\begin{rmk}
Of course, if the action is Zariski locally linearizable then it is \'etale locally linearizable. Celebrated results of Sumihiro imply that if $X$ is a normal variety then any action of $\G_{\mathrm{m},\k}$ is Zariski locally linearizable; see~\cite[\S 0.2]{Richarz20} for references. More recent results of Alper--Hall--Rydh show that \'{e}tale local linearizability is automatic under very mild assumptions on $X$, see again~\cite[\S 0.2]{Richarz20} for references. In any case, for the schemes and actions appearing in the construction of the geometric Satake equivalence, one can check explicitly that the actions are Zariski locally linearizable; see~\S\ref{sss:attractors-Gr} below.
\end{rmk}

The following theorem gathers results from~\cite[Theorem~1.8 and Corollary~1.12]{Richarz16}.

\begin{thm}
\label{thm:attractors-representability}
Assume that the $\G_{\mathrm{m},\k}$-action on $X$ is \'{e}tale locally linearizable. 
Then $X^\circ$, $X^+$ and $X^-$ are schemes over $\k$, and the morphism~\eqref{eqn:morphism-fixed-pts} is a closed immersion. 
Moreover:
\begin{enumerate}
\item
if $X$ is of finite type over $\k$, then so are $X^\circ$, $X^+$ and $X^-$,
\item
the morphisms in~\eqref{eqn:morphisms-att-fixed-pts}
are affine, with geometrically connected fibers, and they induce bijections between sets of connected components of the underlying topological spaces.
\end{enumerate}
\end{thm}

\begin{exm}
\label{ex:attractors}
To get a feeling of what fixed points, attractors and repellers look like, one can consider the following examples.
\begin{enumerate}
\item
\label{ex:attractors-1}
Assume that $X=\Spec(A)$ is an affine $\k$-scheme endowed with an action of $\G_{\mathrm{m},\k}$. The datum of this action is equivalent to the datum of a grading $A = \bigoplus_{n \in \Z} A^n$ such that multiplication is homogeneous (of degree $0$). Any such action is clearly Zariski locally linearizable, and we have the following descriptions:
\begin{itemize}
\item
$X^\circ$ is the closed subscheme of $X$ determined by the ideal of $A$ generated by the subspace $\bigoplus\limits_{n \in \Z \smallsetminus \{0\}} A^n$;
\item
$X^+$ is the closed subscheme of $X$ determined by the ideal of $A$ generated by the subspace $\bigoplus\limits_{n \in \Z_{<0}} A^n$;
\item
$X^-$ is the closed subscheme of $X$ determined by the ideal of $A$ generated by the subspace $\bigoplus\limits_{n \in \Z_{>0}} A^n$.
\end{itemize}
(This---easy---special case is in fact an important ingredient of the proof of Theorem~\ref{thm:attractors-representability}; see~\cite[Lemma~1.9]{Richarz16}.)
\item
\label{ex:attractors-2}
Let $V$ be a $\k$-vector space with a linear action of $\G_{\mathrm{m},\k}$, and consider the corresponding decomposition into eigenspaces: $V=\bigoplus\limits_{j\in \Z} V_j$. For $i \in \Z$ we also set $V_{\geq i} := \bigoplus\limits_{j\geq i} V_j$.
Then the action of $\G_{\mathrm{m},\k}$ on $\P(V)$ is Zariski locally linearizable, and we have 
\[
\P(V)^\circ=\bigsqcup_{\substack{i \in \Z, \\ V_i\neq 0}} \P(V_i) \quad \text{and}\quad  \P(V)^+=\bigsqcup_{\substack{i \in \Z, \\ V_i\neq 0}} \P(V_{\geq i}) \smallsetminus \P(V_{\geq i+1}).
\] 
\end{enumerate}
In case~\ref{ex:attractors-1} the morphisms~\ref{eqn:morphisms-att} are closed immersions, but as illustrated in  case~\ref{ex:attractors-2} this is not the case in general. In fact the latter case illustrates another general phenomenon: if $X$ is a proper $\k$-scheme (with an \'etale locally linearizable action) then the morphisms~\ref{eqn:morphisms-att} induce bijections on the underlying topological spaces, see~\cite[Lemma~1.1.4]{central}. Hence one can think of $X^+$ and $X^-$ as ``the same space as $X$ with a different topology.'' This case is more representative of what happens in the setting of the geometric Satake equivalence.
\end{exm} 

The lemma is an easy exercise in the manipulation of the main results of~\cite{Richarz16}; for details, see~\cite[Lemma~1.1.3]{central}.

\begin{lem}
\label{lem:attractors-closed-subscheme}
Let $X$ be a $\k$-scheme equipped with an \'etale (resp.~Zariski) locally linearizable action of $\G_{\mathrm{m},\k}$, and let  $Y\subset X$ be a $\G_{\mathrm{m},\k}$-stable closed subscheme. 
Then the action of $\G_{\mathrm{m},\k}$ on $Y$ is \'etale (resp.~Zariski) locally linearizable, and moreover we have isomorphisms 
\[
Y^\circ\simeq Y\times_X X^\circ, \quad Y^\pm \simeq Y\times_X X^\pm.
\] 
\end{lem}

Lemma~\ref{lem:attractors-closed-subscheme} allows us to generalize Theorem~\ref{thm:attractors-representability} to ind-schemes, as follows. If $X$ is an ind-scheme over $\k$, an action of $\G_{\mathrm{m},\k}$ on $X$ is the datum of an ``action morphism'' $\G_{\mathrm{m},\k} \times_{\Spec(\k)} X \to X$ which satisfies the obvious axioms. Given such a datum, one can define the functors $X^\circ$, $X^+$ and $X^-$ as in Definition~\ref{def:attractors}.

In practice one often wants to consider actions that ``come from actions on schemes.'' In particular, we will say the action is \'etale (resp.~Zariski) locally linearizable if $X$ admits a presentation $X = \cLim{i} X_i$ such that the action of $\G_{\mathrm{m},\k}$ on $X$ factors through an action on $X_i$ for any $i$, and if moreover the action of $\G_{\mathrm{m},\k}$ on $X_i$ is \'etale (resp.~Zariski) locally linearizable for any $i$. The following result is an easy consequence of Lemma~\ref{lem:attractors-closed-subscheme}; for details, see~\cite[Theorem~1.1.6]{central}.

\begin{thm}
\label{thm:attractors-representability-indsch}
Let $X$ be an ind-scheme over $\k$ equipped with an \'etale locally linearizable of $\G_{\mathrm{m},\k}$, and consider a presentation $X = \cLim{i} X_i$ as in the definition preceding the theorem.
\begin{enumerate}
\item
The functor $X^\circ$ is an ind-scheme over $\k$, which admits a presentation
\[
X^\circ = \cLim{i} \ (X_i)^\circ.
\]
Moreover, the natural morphism $X^\circ \to X$ (similar to~\eqref{eqn:morphism-fixed-pts}) is representable by a closed immersion.
\item
The functor $X^+$ is an ind-scheme over $\k$, which admits a presentation
\[
X^+ = \cLim{i} \ (X_i)^+.
\]
\item
The functor $X^-$ is an ind-scheme over $\k$, which admits a presentation
\[
X^- = \cLim{i} \ (X_i)^-.
\]
\end{enumerate}
\end{thm}

\begin{rmk}
In the setting of Theorem~\ref{thm:attractors-representability-indsch}, the natural morphism $X^+ \to X$ (similar to the left-hand side in~\eqref{eqn:morphisms-att}) is ``representable by schemes'' in the sense that for any scheme $Z$ and any morphism $Z \to X$ the fiber product $X^+ \times_X Z$ is a scheme; see again~\cite[Theorem~1.1.6]{central} for details.
\end{rmk}

\subsubsection{Braden's Theorem} 
\label{sss:Braden}

By a ``ring of coefficients'' we mean a ring of one of the following forms:
\begin{itemize}
\item
a finite field of characteristic invertible in $\k$, or an algebraic closure of such a field;
\item
the ring of integers in a finite extension of $\Q_p$, where $p$ is a prime number invertible in $\k$;
\item
a finite extension or an algebraic closure of $\Q_p$, where $p$ is a prime number invertible in $\k$.
\end{itemize}
If $\Lambda$ is such a ring, for any $\k$-scheme of finite type we will denote by $\Dbc(Y,\Lambda)$ the bounded derived category\footnote{We make the usual abuse of terminology in this setting: unless $\Lambda$ is a finite field, this category is \emph{not} the bounded derived category of the abelian category of \'etale $\Lambda$-sheaves, but something more complicated constructed using an appropriate limit procedure.} of \'etale $\Lambda$-sheaves on $Y$.

Fix now $X$ a scheme of finite type over $\k$, endowed with an \'{e}tale locally linearizable action of $\G_{\mathrm{m},\k}$, and a ring of coefficients $\Lambda$.
Recall the morphisms $q^\pm$ and $p^\pm$ considered in~\eqref{eqn:morphisms-att-fixed-pts}--\eqref{eqn:morphisms-att}.
We define the \emph{hyperbolic localization functors} as the following functors from $\Dbc(X,\Lambda)$ to $\Dbc(X^\circ,\Lambda)$
\[
\mathrm{HL}^+_X= (q^+)_! (p^{+})^* , \quad \mathrm{HL}^-_X= (q^-)_* (p^{-})^! .
\]
As explained in~\cite[Construction~2.2]{Richarz16}, there exists a canonical morphism of functors
\begin{equation}
\label{eqn:morph-Braden}
\mathrm{HL}^-_X\rightarrow \mathrm{HL}^+_X. 
\end{equation}

Let $a,p \colon \G_{\mathrm{m}, \k} \times_{\Spec(\k)} X\rightarrow X$ be the action and projection maps, respectively.
We will denote by 
\[
\Dbc(X,\Lambda)^{\G_{\mathrm{m},\k}\mhyphen\mathrm{mon}}
\]
the full triangulated subcategory of $\Dbc(X,\Lambda)$ generated by the complexes $\scF$ such that there exists an isomorphism $a^*\scF\simto p^*\scF$. (This condition should be seen as a ``weak equivariance'' condition.)

The following theorem is proved in~\cite[Theorem~2.6]{Richarz16}. Earlier variants are due to Braden and Drinfeld--Gaitsgory; see~\cite{Richarz16} for precise references.

\begin{thm}
\label{thm:braden}
For any $\scF\in \Dbc(X,\Lambda)^{\G_{\mathrm{m},\k}\mhyphen\mathrm{mon}}$, the morphism 
\[
\mathrm{HL}^-_X(\scF) \rightarrow \mathrm{HL}^+_X(\scF)
\]
induced by~\eqref{eqn:morph-Braden}
is an isomorphism. 
\end{thm}

\begin{rmk}
\phantomsection
\label{rmk:HL}
\begin{enumerate}
\item 
The morphism~\eqref{eqn:morphism-fixed-pts} factors natural through morphisms
\[
i^+ \colon X^\circ \to X^+ \quad \text{and} \quad i^- : X^\circ \to X^-
\]
(given by ``constant maps''). As explained in~\cite[(2.4)]{Richarz16}, there exists canonical morphisms of functors from $\Dbc(X^-, \Lambda)$ to $\Dbc(X^\circ, \Lambda)$, resp.~from $\Dbc(X^+, \Lambda)$ to $\Dbc(X^\circ, \Lambda)$,
\[
(q^-)_* \to (i^-)^*, \quad \text{resp.} \quad (i^+)^! \to (q^+)_!.
\]
There exists natural actions of $\G_{\mathrm{m},\k}$ on $X^-$ and $X^+$ (see the comments preceding \ref{rmk:att-monomorphism-separated}), so that as above one can consider the subcategories $\Dbc(X^-,\Lambda)^{\G_{\mathrm{m},\k}\mhyphen\mathrm{mon}}$ and $\Dbc(X^+,\Lambda)^{\G_{\mathrm{m},\k}\mhyphen\mathrm{mon}}$. If $\scF$ belongs to $\Dbc(X^-,\Lambda)^{\G_{\mathrm{m},\k}\mhyphen\mathrm{mon}}$, resp.~if $\scF'$ belongs to $\Dbc(X^+,\Lambda)^{\G_{\mathrm{m},\k}\mhyphen\mathrm{mon}}$, then the morphism above induces an isomorphism
\[
(q^-)_* \scF \simto (i^-)^* \scF, \quad \text{resp.} \quad (i^+)^! \scF' \simto (q^+)_! \scF'.
\]
In fact this claim can be checked \'etale locally, so that using~\cite[Lemma~1.11]{Richarz16} one can assume that $X$ is affine.
This claim is proved in~\cite[Lemma~2.19]{Richarz16} in case $X$ is a vector space with a linear action. Using~\cite[Lemmas~2.21,~2.23,~2.25]{Richarz16}, we deduce that its holds for any affine $X$, as desired.

Using this claim, we obtain that
for any $\scF\in \Dbc(X,\Lambda)^{\G_{\mathrm{m},\k}\mhyphen\mathrm{mon}}$ we have canonical isomorphisms 
\begin{equation}
\label{eqn:HL-push-pull}
\mathrm{HL}^-_X(\scF) \simeq (i^-)^* (p^{-})^! \scF \quad \text{and} \quad 
\mathrm{HL}^+_X(\scF) \simeq (i^+)^! (p^{+})^* \scF.
\end{equation}
Hence one can consider that the hyperbolic localization functors are obtained by ``taking a $!$-pullback functor in some directions, followed by a $*$-pullback functor in other directions.''
\item 
\label{it:HL-commutes}
In practice, one can often prove that appropriate functors commute (in the appropriate sense) with the hyperbolic localization functor by using the flexibility to pass between $!$- and $*$-functors using the isomorphism in Theorem~\ref{thm:braden}, or push and pull functors using the isomorphisms~\eqref{eqn:HL-push-pull}. For concrete applications of this idea, see~\cite[Proposition~3.1 and Theorem~3.3]{Richarz16}.
\end{enumerate}
\end{rmk}

\section{The geometric Satake equivalence} 
\label{sec:lecture-2}

In this section, we introduce the geometric Satake equivalence. What we understand of the history of this construction is discussed in full detail in the introduction of~\cite{central}, hence will not be repeated here. We will only cite the names of the main contributors in this story, in chronological order, namely Lusztig, Ginzburg, Be{\u\i}linson--Drinfeld and Mirkovi{\'c}--Vilonen. The first complete proof of this equivalence, which is also the only one known in the generality that we require, is to be found in~\cite{MV07}. See the introduction for references to more detailed surveys of this proof.


\subsection{More on the geometry of affine Grassmannians}

\subsubsection{Sheaves on affine Grassmannians} 
\label{ss:sheaves-Gr}

Let $G$ be a smooth affine group scheme of finite type over $\D_\k$. 
Recall (see Proposition~\ref{prop:loop-functors}) that 
\[
\L^+ G= \varprojlim_{i} \L^+_i G
\]
is an affine group scheme over $\k$, with each $\L^+_i G$ smooth affine of finite type. It is also known that for any $i \geq 1$ the morphism $\L^+_{i+1} G \to \L^+_i G$ is surjective, and that its kernel is a vector group, i.e.~the $\k$-algebraic group associated with a finite-dimensional vector space over $\k$, hence in particular a smooth connected unipotent group. (See e.g.~\cite[\S VI.1]{fargues-scholze} for a statement of this form.)


Consider also the affine Grassmannian $\GGr_G$, which is a separated ind-scheme of ind-finite type, see Theorem~\ref{thm:GrG-indscheme}.
We have an action of $\L G$, hence of $\Loop^+ G$, on $\GGr_G$ by twisting the section $\alpha$. 
In fact, a look at the proof of this theorem shows that we have a presentation
\begin{equation}
\label{eqn:presentation-GrG}
\GGr_G=\cLim{i} \GGr_{G,i}
\end{equation}
such that each $\GGr_{G,i}$ is separated of finite type, stable under the action of $\L^+ G$, and moreover the action factors through an action of $\L^+_{n_i} G$ for some $n_i\geq 0$. 
Then if $\Lambda$ is a ring of coefficients (see~\S\ref{sss:Braden}), for any $m \geq n_i$ we can consider the constructible $\L^+_m G$-equivariant derived category $\Db_{\L^+_{m} G}(\GGr_{G,i}, \Lambda)$ of \'etale $\Lambda$-sheaves on $\GGr_{G,i}$.
(This category can be constructed either using the Bernstein--Lunts construction, or using the theory of \'etale sheaves on stacks, applies to the algebraic stack of finite type $\L^+_{m} G \backslash \GGr_{G,i}$.) Using the fact that the quotient morphisms $\L^+_{j+1} G \to \L^+_j G$ have connected unipotent kernels, one sees that this category does not depend on the choice of $m$ up to canonical equivalence, hence can be denoted
\[
\Db_{\L^+ G}(\GGr_{G,i}, \Lambda).
\]
If $i\leq j$ we have a fully-faithful functor
\[
\Db_{\L^+ G}(\GGr_{G,i}, \Lambda)\rightarrow \Db_{\L^+ G}(\GGr_{G,j}, \Lambda)
\]
given by pushforward along the closed immersion $\GGr_{G,i} \to \GGr_{G,j}$.
Therefore, we can consider the \emph{derived Satake category}, given as the colimit
\[
\Db_{\L^+ G}(\GGr_{G}, \Lambda):= \cLim{i\geq 0} \Db_{\L^+ G}(\GGr_{G,i}, \Lambda).
\]
Using~\eqref{eqn:Hom-ind-sch} one can check that this category
does not depend (up to canonical equivalence) on the choice of the presentation~\eqref{eqn:presentation-GrG}
satisfying the properties stated above. Below we will also sometimes consider the category
\[
\Dbc(\GGr_{G}, \Lambda):= \cLim{i\geq 0} \Dbc(\GGr_{G,i}, \Lambda).
\]

In practice, to avoid trivial but cumbersome discussions, one often works with this category by ``pretending that $\Gr_G$ is a scheme,'' i.e.~writing formulas that do not actually make sense as is, but would make sense if all the (ind-)schemes under consideration were schemes of finite type, and can be easily modified to make actual sense. A similar comment applies to other ind-schemes to be introduced below.

For any $i \geq 0$ the category $\Dbc(\GGr_{G,i}, \Lambda)$ has a canonical perverse t-structure, and there exists a unique t-structure on $\Db_{\L^+ G}(\GGr_{G,i}, \Lambda)$ such that the natural forgetful functor
\[
\Db_{\L^+ G}(\GGr_{G,i}, \Lambda) \to \Dbc(\GGr_{G,i}, \Lambda)
\]
is t-exact.\footnote{Often, the definition of the perverse t-structure on sheaves on a stacks is normalized in such a way that for a quotient stack $X/H$ the pullback functor $\Dbc(X/H, \Lambda) \to \Dbc(X,\Lambda)$ is t-exact \emph{up to a shift}. In the present context it is more convenient to add this shift \emph{in the definition of the t-structure}, so that the pullback functor becomes t-exact.}
Then there exists a unique t-structure on $\Db_{\L^+ G}(\GGr_{G}, \Lambda)$ such that each natural functor
\[
\Db_{\L^+ G}(\GGr_{G,i}, \Lambda) \to \Db_{\L^+ G}(\GGr_{G}, \Lambda)
\]
is t-exact. Again, this t-structure does not depend on the choice of presentation~\eqref{eqn:presentation-GrG}. The heart of this t-structure will be denoted $\Perv_{\L^+ G}(\GGr_G,\Lambda)$, and the associated cohomology functors will be denoted $(\pH^n : n \in \Z)$. Similar comments apply in the category $\Dbc(\GGr_{G}, \Lambda)$.

\begin{rmk}
\label{rmk:etale-const}
In these notes we work with \'etale sheaves, because this is the setting we will require for our main application later. However the geometric Satake equivalence can also be stated (and proved) in the setting where $\k=\mathbb{C}$, working with sheaves of $\Lambda$-modules on the associated analytic space. This is the setting considered explicitly in~\cite{MV07,BR18}. The two settings can be treated with essentially identical methods.
\end{rmk}

\subsubsection{Attractors} 
\label{sss:attractors-Gr}

From now on, we assume that $\k$ is algebraically closed, and
fix a connected reductive algebraic group $G$ over $\k$. For a moment, we will only consider the constructions above applied to the group $G \times_{\Spec(\k)} \D_\k$.

We also fix opposite Borel subgroups $B^+$, $B^-$, and denote by $T$ their intersection (a maximal torus in $G$).
Let $W=N_G(T)/T$ be the Weyl group of $(G,T)$. 
We will also consider the character and cocharacter lattices
\[
\X=X^*(T), \quad \X^\vee=X_*(T),
\]
and denote the root and coroot systems by $R\subset \X$ and $R^\vee\subset \X^\vee$ respectively. The $T$-weights in the Lie algebra of $B^+$ determine a system of positive roots $R_+ \subset R$, and we denote by
$\X^\vee_+\subset \X^\vee$ the corresponding set of dominant coweights. 
There is a partial order $\leq$ on $\X^\vee$ such that $\lambda \leq \mu$ if and only if $\mu-\lambda\in \Z_{\geq 0}R^\vee$. 
W set $\rho=\frac{1}{2}\sum_{\alpha\in R_+} \alpha \in \Q \otimes_{\Z} \X$. 

Let us fix a strictly dominant coweight, i.e.~an element
$\nu\in \X^\vee$ such that $\langle \nu ,\alpha \rangle >0$ for any $\alpha\in R_+$. 
Then we have a $\G_{\mathrm{m},\k}$-action on $G$ given by for $g\in G$ and $t\in \G_{\mathrm{m},\k}$ by $t \cdot g=\nu(t)g\nu(t)^{-1}$. 
It is a classical fact that for this action the attractor, repeller and fixed points can be described as follows: 
\[
G^+=B^+, \quad G^-=B^-, \quad G^\circ=T;
\] 
see e.g.~\cite[\S 2.1]{CGP10}.
We also have an action of $\G_{\mathrm{m},\k}$ on the ind-scheme $\Gr_G$ obtained by restriction of the action of $\L G$ via the composition $\G_{\mathrm{m},\k} \xrightarrow{\nu} G \hookrightarrow \Loop^+G$, where the second morphism is given by taking constant loops.

The following statement is implicit in~\cite{MV07}, and was proved in this form by Haines--Richarz. For a description of this proof, and precise references, see~\cite[\S 1.2.2]{central}.

\begin{thm} 
\phantomsection
\label{thm:attractors-Gr}
\begin{enumerate}
\item The action of $\G_{\mathrm{m},\k}$ on $\Gr_G$ is Zariski locally linearizable. 
\item We have canonical identifications 
\[
\Gr_{B^+}\simeq (\Gr_G)^+, \quad \Gr_{B^-}\simeq (\Gr_G)^- \quad \text{and} \quad \Gr_T\simeq (\Gr_G)^\circ.
\]
\end{enumerate}
\end{thm}

In particular, this theorem and the comments in Example~\ref{ex:attractors} show that the natural maps $\Gr_{B^+} \to \Gr_G$ and $\Gr_{B^-} \to \Gr_G$ are bijections (but not isomorphisms, and not even homeomorphisms).


\subsubsection{Schubert cells}
\label{sss:schubert-cells}

For each $\lambda \in \X^\vee$ we consider the composition
\[
\Loop \G_{\mathrm{m},\k} \xrightarrow{\Loop \lambda} \Loop G \rightarrow \Gr_G
\]
where the second map is the quotient morphism (see Theorem~\ref{thm:Gr-quotient}).
The element $z\in \Loop \G_{\mathrm{m},\k}(\k)=\k \pa{z}^\times$ maps to a $\k$-point of $\Gr_G$, denoted $\Loop_\lambda \in \Gr_G(\k)$. 
For $\mu\in \X^\vee_+$ we will denote by
\[
\Gr_{G,\leq \mu}
\]
the closure of the $\Loop^+ G$-orbit of $\Loop_\mu$, endowed with the reduced scheme structure. (This definition makes sense because $\Gr_G$ is a colimit of schemes on which the action factors through a quotient of finite type.)
This scheme is a projective $\k$-variety. 

It is a standard fact that for $\mu, \mu'\in \X^\vee_+$, we have 
\[
\Gr_{G,\leq \mu'}\subset \Gr_{G,\leq \mu} \quad \text{if and only if} \quad \mu\leq \mu'.
\]
For $\mu \in \X^\vee_+$ we set 
\[
\Gr_{G,\mu}:=\Gr_{G,\leq \mu}\backslash \bigcup_{\substack{\mu' \in \X^\vee_+ \\ \mu'< \mu}} \Gr_{G,\leq \mu'}.
\]
Then $\Gr_{G,\mu}$ is an open subscheme in $\Gr_{G,\leq \mu}$, hence a quasi-projective scheme over $\k$.
It is known that $\Gr_{G,\mu}$ is smooth, of dimension $\langle \mu, 2\rho\rangle$, and that its $\k$-points consist of the orbit of $\Loop_\mu$ under the action of $\Loop^+ G (\k)$.
Moreover the stabilizer of $\Loop_\mu$ in $\Loop^+G$ is connected. 
In particular, for the underlying topological spaces we have a stratification
\[
|\Gr_G|= \bigsqcup_{\lambda \in \X^\vee_+} |\Gr_{G,\lambda}|.
\]
The decomposition of $\Gr_G$ into its connected components is of the form
\[
\Gr_G= \bigsqcup_{c\in \X^\vee/\Z R^\vee} \Gr^c_G,
\]
where $\Gr_{G,\mu}\subset \Gr^c_G$ if and only if $\mu\in c$. 

\subsubsection{Semi-infinite orbits} 
\label{sss:si-orbits}

Consider the morphisms 
\begin{equation}
\label{eqn:morphisms-Gr-B-T}
\Gr_{B^+}\rightarrow \Gr_T, \quad \Gr_{B^-}\rightarrow \Gr_T
\end{equation}
induced by the natural morphisms $B^+ \rightarrow T$, $B^- \to T$ (given by the quotient by the unipotent radical). 
It is a standard fact that $|\Gr_T|$ is discrete, with underlying set $\X^\vee$.
It can be easily deduced from Theorem~\ref{thm:attractors-representability} that the maps in~\eqref{eqn:morphisms-Gr-B-T} induce bijections between the corresponding sets of connected components; see~\cite[\S 1.2.2.3]{central} for details. In this way, we obtain a bijection between the set of connected components of $| \Gr_{B^+} |$, resp.~$| \Gr_{B^-} |$, and $\X^\vee$.
For any $\lambda\in \X^\vee$, we will denote by
\[
\rS_\lambda, \quad \text{resp.} \quad \rT_\lambda,
\]
the connected component\footnote{See~\cite[\S 1.1.1.3]{central} for comments on the ind-scheme structure on connected components of ind-schemes.} of $\Gr_{B^+}$, resp.~$\Gr_{B^-}$, corresponding to $\lambda$. 
Then we have morphisms 
\begin{equation}
\label{eqn:embedding-S-T}
s_\lambda\colon \rS_\lambda \rightarrow \Gr_G \quad \text{and}\quad t_\lambda\colon \rT_\lambda \rightarrow \Gr_G
\end{equation}
which are representable by locally closed immersions. 
We also have 
\[
\rS_\lambda(\k)=B^+\big(\k\pa{z} \big) \cdot \Loop_\lambda= U^+\big(\k\pa{z} \big) \cdot \Loop_\lambda, \quad
\rT_\lambda(\k)=B^-\big(\k\pa{z} \big) \cdot \Loop_\lambda= U^-\big(\k\pa{z} \big) \cdot \Loop_\lambda,
\]
where $U^+$ is the unipotent radical of $B^+$ and $U^-$ is the unipotent radical of $B^-$.

\subsubsection{Dimension estimate} 

For $\lambda\in \X^\vee$, we denote by $\mathsf{dom}(\lambda)$ the unique dominant $W$-conjugate of $\lambda$.
For any $\mu\in \X^\vee_+$, we set $\X^\vee_\mu :=\{\lambda\in \X^\vee \mid \mathsf{dom}(\lambda) \leq \mu\}$. The following statement is a reformulation of a crucial result in~\cite{MV07}; see~\cite[Theorem~1.5.2]{BR18},~\cite[Remark~1.2.9]{central} and~\cite[Talk~8]{workshop}.

\begin{thm}
\label{thm:dimension-estimate}
For any $\lambda\in \X^\vee$ and $\mu\in \X^\vee_+$, the fiber product 
\[
\rS_\lambda\times_{\Gr_G} \Gr_{G,\leq \mu}
\]
is a connected affine scheme of finite type, which is nonempty if and only if $\lambda\in \X^\vee_\mu$. 
In this case, each of its irreducible components has dimension $\langle \rho, \mu+\lambda\rangle$. 
\end{thm}

\subsection{The geometric Satake equivalence}

\subsubsection{Dual group}

The root datum of the pair $(G,T)$ is the quadruple $(\X,\X^\vee, R, R^\vee)$, together with the natural bijection between roots and coroots and the natural perfect pairing between $\X$ and $\X^\vee$. One can obtain another root datum $(\X^\vee,\X, R^\vee, R)$ by switching the roles of weights/roots and coweights/coroots.
By the work of Chevalley, Demazure and Grothendieck, there exists a unique split reductive group scheme $G^\vee_\Z$ over $\Spec(\Z)$ with a split maximal torus $T^\vee_\Z$, such that for any algebraically closed field $\F$, the pair $(G^\vee_\F,T^\vee_\F)$ obtained by base change has root datum $(\X^\vee,\X, R^\vee, R)$; see~\cite[\S 1.14.1]{BR18} for comments and references. For any commutative ring $\Lambda$ we set 
\[
G^\vee_\Lambda=\Spec(\Lambda) \times_{\Spec(\Z)} G^\vee_\Z,
\]
seen as a (reductive) group scheme over $\Spec(\Lambda)$.

\subsubsection{Convolution}
\label{sss:convolution}

Let us consider the ``convolution affine Grassmannian" 
\[
\Conv_G :=(\Loop G\times \Gr_G \big/ \Loop^+G)_{\mathrm{\acute{e}t}},
\]
where $\Loop^+G$ acts on $\Loop G\times \Gr_G$ via the formula $g \cdot (h,x)=(hg^{-1},gx)$ for $g \in \Loop^+G$, $h \in \Loop G$ and $x \in \Gr_G$. Here the automorphism of $\Loop G\times \Gr_G$ given by $(h,x) \mapsto (h, h \cdot x)$ identifies this action with the action induced by multiplication on the right on the left-hand factor, so that this \'etale quotient is an ind-scheme (isomorphic to $\Gr_G \times \Gr_G$, but which should not be confused with it). It is not difficult to describe a moduli problem (in terms of torsors) that this ind-scheme represents; see e.g.~\cite[\S 1.7.2]{BR18}. If $\Lambda$ is a ring of coefficients, one defines the $\Loop^+ G$-equivariant derived category $\Db_{\Loop^+ G}(\Conv_G,\Lambda)$ as in the case of the affine Grassmannian in~\S\ref{ss:sheaves-Gr}.
 
The morphism $\Loop G\times \Gr_G \to \Gr_G$ given by $(g,x) \mapsto g \cdot x$ factors through a morphism
\[
m\colon \Conv_G \rightarrow \Gr_G.
\] 
Consider also the natural projection morphism $p\colon \Loop G \rightarrow \Gr_G$. 
Given complexes $\scF,\scG\in \Db_{\Loop^+ G}(\Gr_G,\Lambda)$, there exists a unique complex 
\[
\scF \tboxtimes \scG \quad \in \Db_{\Loop^+ G}(\Conv_G,\Lambda)
\]
whose pullback to $\Loop G \times \Gr_G$ is $p^*\scF \lboxtimes_\Lambda \scG$ (as an equivariant complex). 
We set
\[
\scF \star \scG:= m_*(\scF \tboxtimes \scG).
\]
It is a standard fact that the bifunctor $\star$ defines a monoidal structure on the category $\Db_{\Loop^+ G}(\Gr_G,\Lambda)$. 

The following property is crucial; for a discussion of the proof, see~\cite[\S 1.10.3]{BR18}.

\begin{prop}
If $\scF,\scG\in \Perv_{\Loop^+ G}(\Gr_G,\Lambda)$, then we have 
\[
\pH^n(\scF \star \scG)=0 \quad \text{if $n>0$.}
\] 
\end{prop}

An easy consequence of this proposition is that
the bifunctor $\star^0$ on $\Perv_{\Loop^+G}(\Gr_G,\Lambda)$ defined by 
\[
\scF \star^0 \scG= \pH^0(\scF \star \scG)
\]
admits canonical unitality and associativity constraints which define a monoidal structure on the abelian category $\Perv_{\Loop^+G}(\Gr_G,\Lambda)$. (The monoidal product is exact in case $\Lambda$ is a field, but not otherwise.)

\subsubsection{Statement}

We can finally state the geometric Satake equivalence.
Fix a ring of coefficients $\Lambda$, and
denote by $\Rep(G^\vee_\Lambda)$ the category of $G^\vee_\Lambda$-representations (i.e., $\mathscr{O}(G^\vee_\Lambda)$-comodules) that are finitely generated as $\Lambda$-modules. We will also denote by
\[
\mathrm{For}_{G^\vee_\Lambda} \colon \Rep(G^\vee_\Lambda) \to \Lambda\mhyphen\mathrm{Mod}
\]
the functor of forgetting the action.

\begin{thm}
\label{thm:Satake}
There exists an equivalence of monoidal categories 
\[
\mathrm{Sat}\colon \big(\Perv_{\Loop^+G}(\Gr_G,\Lambda), \star^0\big) \simto \big(\Rep(G^\vee_\Lambda), \otimes_\Lambda \big)
\]
together with an isomorphism of functors
\[
\mathrm{For}_{G^\vee_\Lambda}\circ \mathrm{Sat}\simeq \coH^\bullet(\Gr_G,-).
\]
\end{thm}

\subsubsection{Commutativity constraint (fusion product)} 

Note that the monoidal category 
\[
\big(\Rep(G^\vee_\Lambda), \otimes_\Lambda \big)
\]
has a commutativity constraint. 
Hence in view of Theorem~\ref{thm:Satake}, the monoidal category 
\[
\big(\Perv_{\Loop^+G}(\Gr_G,\Lambda), \star^0 \big)
\]
should also have a commutativity constraint. 
In fact, the prior construction of this structure is a crucial step in the \emph{proof} of Theorem~\ref{thm:Satake}. 

To construct the desired commutativity constraint, one needs to work over the curve 
\[
C=\A^1_\k=\Spec (\k[x]).
\]
For $R\in \Alg_\k$, we set $C_R := C\times_{\Spec(\k)} \Spec(R)$. 
For any $y\in C(R)$ (resp. $y_1,y_2\in C(R)$) we denote by $\widehat{\Gamma}_y$ (resp. $\widehat{\Gamma}_{y_1,y_2}$) the spectrum of the completion of $R[x]$ with respect to the ideal defining the graph $\Gamma_y$ of $y$ (resp.~the union $\Gamma_{y_1}\cup \Gamma_{y_2}$). 
Variants of Theorem~\ref{thm:GrG-indscheme} 
imply that there exist ind-schemes $\Gr_{G,C}$ and $\mathrm{Fus}_G$ such that, for $R \in \Alg_\k$,
\begin{itemize}
\item
$\Gr_{G,C}(R)$ is given by the set of isomorphism classes of triples $(y,\sE,\beta)$ where $y\in C(R)$, $\sE$ is a $(G \times_{\Spec(\k)} \widehat{\Gamma}_y)$-torsor, and $\beta$ is a trivialization of this torsor on $\widehat{\Gamma}_y \backslash \Gamma_y$;
\item
$\mathrm{Fus}_{G}(R)$ is given by the set of isomorphism classes of quadruples $(y_1,y_2,\sE,\beta)$ where $y_1, y_2 \in C(R)$, $\sE$ is a $(G \times_{\Spec(\k)} \widehat{\Gamma}_{y_1,y_2})$-torsor, and $\beta$ is a trivialization of this torsor on $\widehat{\Gamma}_{y_1,y_2} \backslash (\Gamma_{y_1}\cup \Gamma_{y_2})$.
\end{itemize}
These ind-schemes are equipped with canonical morphisms $\Gr_{G,C}\rightarrow C$ and $\mathrm{Fus}_G\rightarrow C^2$.


The following statement is standard; see~\cite[\S 1.7.3]{BR18} for a discussion. Here we denote by $\Delta C$ the diagonal copy of $C$ in $C^2$.

\begin{prop}
\begin{enumerate}
	\item There exists a canonical isomorphism $\Gr_{G,C} \simeq \Gr_G \times C$. 
	\item We have isomorphisms 
\[
\mathrm{Fus}_G \times_{C^2} \Delta C \simeq \Gr_{G,C}, \quad
\mathrm{Fus}_G \times_{C^2} (C^2\backslash \Delta C) \simeq (\Gr_{G,C}\times \Gr_{G,C}) \times_{C^2} (C^2 \smallsetminus \Delta C).
\]
\end{enumerate}
\end{prop}

In view of this proposition we have a canonical morphism representable by an open immersion
\[
j\colon \Gr_{G}\times \Gr_{G}\times (C^2\smallsetminus \Delta C) \hookrightarrow \mathrm{Fus}_G,
\]
and a canonical morphism representable by a closed immersion
\[
i\colon \Gr_{G,C} \hookrightarrow \mathrm{Fus}_G.
\]
We will consider the middle extension functor $j_{!*}$ associated with $j$, and the pullback functor $i^*$ associated with $i$.

For the following statement, we refer to~\cite[Lemma~1.7.10]{BR18}.

\begin{prop}
For $\scF,\scG\in \Perv_{\Loop^+ G}(\Gr_G, \Lambda)$ there is a canonical isomorphism 
\[
i^*j_{!*}\big(
\pH^0(\scF \lboxtimes_{\Lambda} \scG) \lboxtimes_{\Lambda} \underline{\Lambda}_{C^2 \smallsetminus \Delta C}[2] \big) 
\simeq (\scF \star^0 \scG) \lboxtimes_{\Lambda} \underline{\Lambda}_{\Delta_C}.
\]
\end{prop} 

The description of the convolution product obtained from this proposition is the crucial ingredient in the construction of the commutativity constraint for the product $\star^0$, which takes advantage of the possibility to switch the factors in the product $\Gr_{G}\times \Gr_{G}\times (C^2\smallsetminus \Delta C)$.

\subsubsection{Weight functors}
\label{sss:weight-functors}

Recall the morphisms~\eqref{eqn:embedding-S-T}. Given $\lambda \in \X^\vee$ and an object $\scF \in \Perv_{\Loop^+ G}(\Gr_G, \Lambda)$ we consider the graded $\Lambda$-modules
\[
\coHc^\bullet(\rS_\lambda ,s_\lambda^*\scF)\quad 
\text{and}\quad \coH^\bullet_{\rT_\lambda}(\scF):=\coH^\bullet(\rT_\lambda, t_\lambda^!\scF).
\]
In view of Theorem~\ref{thm:attractors-Gr}, the direct sums over all $\lambda$ of these functors are hyperbolic localization functors in the sense of~\S\ref{sss:Braden}.

For the proof of the following statement, we refer to~\cite[\S\S 1.5.3--1.5.4]{BR18}.

\begin{thm}
\phantomsection
\label{thm:weight-functors}
\begin{enumerate}
\item 
\label{it:weight-functors-1}
For any $\lambda\in \X^\vee$ and $\scF\in \Perv_{\Loop^+G}(\Gr_G,\Lambda)$, there exists a canonical isomorphism of graded $\Lambda$-modules
\[
\coHc^\bullet(\rS_\lambda, s_\lambda^*\scF) \simeq \coH^\bullet_{\rT_\lambda}(\scF).
\] 
Moreover, these graded modules vanish in all degrees except possibly $\langle\lambda,2\rho\rangle$. 
\item 
\label{it:weight-functors-2}
For $\scF\in \Perv_{\Loop^+G}(\Gr_G,\Lambda)$, there exists canonical isomorphisms 
\[
\coH^\bullet(\Gr_G, \scF) \simeq \bigoplus_{\lambda\in \X^\vee} \coHc^{\langle\lambda, 2\rho\rangle}(\rS_\lambda, \scF) \simeq \bigoplus_{\lambda\in \X^\vee} \coH^{\langle\lambda, 2\rho\rangle}_{\rT_\lambda}(\scF). 
\]
\end{enumerate}
\end{thm}

In this theorem, the isomorphism in~\eqref{it:weight-functors-1}
follows from Theorem~\ref{thm:braden}. For the vanishing statement, one sees that the computation of the
dimension of $\rS_\lambda \cap \Gr_{G, \leq \mu}$ in Theorem~\ref{thm:dimension-estimate} implies that $\coHc^\bullet(\rS_\lambda, s_\lambda^*\scF)$ vanishes in degrees larger than $\langle \lambda, 2\rho \rangle$, while the computation of the dimension of $\rT_\lambda \cap \Gr_{G, \leq \mu}$ in the same statement implies that $\coH^\bullet_{\rT_\lambda}(\scF)$ vanishes in degrees less that $\langle \lambda, 2\rho \rangle$.
The isomorphisms in~\eqref{it:weight-functors-2} are obtained as a consequence.

Using standard ``long exact sequence'' arguments, the vanishing statement in~\eqref{it:weight-functors-1} implies that the functors
\[
\coHc^{\langle\lambda,2\rho\rangle}(\rS_\lambda, s_\lambda^* (-)) \simeq \coH^{\langle\lambda,2\rho\rangle}_{\rT_\lambda}(-)
\]
are exact. (These functors are called ``weight functors.'')

\subsubsection{Outline of the proof} 

The proof of Theorem~\ref{thm:Satake} is based on ideas from the \emph{tannakian formalism}: one constructs structures on the category $\Perv_{\Loop^+ G}(\Gr_G, \Lambda)$ that will eventually ``force it'' to be the category of representations of a group scheme, and then one shows that this group scheme must be $G^\vee_\Lambda$.

There are some general ``tannakian reconstruction'' results asserting that monoidal categories possessing a number of structures are necessarily categories of representations of a group scheme; see~\cite[\S 1.2]{BR18} for a discussion. These results play a role in the proof of Theorem~\ref{thm:Satake}, but they usually require the category to be linear over a field. The possibility to treat ``integral'' coefficients is crucial for the proof of Theorem~\ref{thm:Satake} (see Remark~\ref{rmk:remarks-proof-Satake} below), so that the ``tannakian reconstruction'' that is performed for the proof of that theorem is done ``by hand,'' exploiting special features of the situation, rather than by an application of a general statement. Namely, if $Z\subset \Gr_G$ is a closed subscheme which is a finite union of $\Loop^+G$-orbits, one checks that the functor 
\[
\Perv_{\Loop^+ G}(Z,\Lambda)\rightarrow \Perv_{\Loop^+ G}(\Gr_G,\Lambda) \xrightarrow{\coH^{\langle\lambda, 2\rho\rangle}_{\rT_\lambda}(-)} \Lambda\mhyphen\mathrm{Mod}
\]
(where the first map is given by pushforward)
is representable by an object $\scP_Z(\lambda)$, see~\cite[\S 1.12.1]{BR18}.
Then one considers 
\[
\scP_Z=\bigoplus_{\substack{\nu\in \X^\vee,\\ Z\cap \rT_\nu\neq \varnothing}} \scP_Z(\nu),
\]
which is shown to be a projective generator of the category $\Perv_{\Loop^+ G}(Z,\Lambda)$. 
Once this is known, considering the algebra $\mathrm{A}_Z=\End(\scP_Z)^\op$ we have an equivalence of abelian categories
\[
\Perv_{\Loop^+ G}(Z,\Lambda) \simto \mathrm{A}_Z\mhyphen\mathrm{Mof}
\]
where the right-hand side is the category of finitely generated $\mathrm{A}_Z$-modules; see~\cite[\S 1.13.1]{BR18}. One checks also that $\mathrm{A}_Z$ is free of finite rank over $\Lambda$, so that the dual $\Lambda$-module $\mathrm{B}_Z:=\Hom_\Lambda(\mathrm{A}_Z,\Lambda)$ has a canonical structure of $\Lambda$-coalgebra, and we have a canonical equivalence of categories
\[
\mathrm{A}_Z\mhyphen\mathrm{Mof} \simeq \mathrm{Comof}\mhyphen\mathrm{B}_Z
\]
where the right-hand side is the category of right $\mathrm{B}_Z$ which are finitely generated over $\Lambda$.


If $Y$ is another closed subscheme as above which contains $Z$, then there is a canonical injective morphism $\mathrm{B}_Z \hookrightarrow \mathrm{B}_Y$ such that the pushforward functor
\[
\Perv_{\Loop^+ G}(Z,\Lambda) \to \Perv_{\Loop^+ G}(Y,\Lambda)
\]
corresponds under the equivalences above to the functor
\[
\mathrm{Comof}\mhyphen\mathrm{B}_Z \to \mathrm{Comof}\mhyphen\mathrm{B}_Y
\]
induced by this morphism.
One can therefore consider the coalgebra 
\[
\mathrm{B}=\cLim{Z} \mathrm{B}_Z,
\]
and we have an equivalence of abelian categories
\[
\Perv_{\Loop^+ G}(\Gr_G,\Lambda) \simto \mathrm{Comof}\mhyphen\mathrm{B}.
\]
From the monoidal structure on $\Perv_{\Loop^+ G}(\Gr_G,\Lambda)$ and its commutativity constraint, we obtain a commutative algebra structure on $\mathrm{B}$ which makes it a Hopf algebra (so that $\Spec(\mathrm{B})$ is a group scheme over $\Spec(\Lambda)$) and equivalences of monoidal categories 
\[
\big(\Perv_{\Loop^+ G}(\Gr_G,\Lambda), \star^0\big) \simto \big(\mathrm{Comof}\mhyphen\mathrm{B}, \otimes_\Lambda \big) \simto \big( \Rep(\Spec(\mathrm{B})), \otimes_\Lambda \big).
\]
Finally, one checks that the formation of $\mathrm{B}$ is compatible in the obvious way with extension of scalars, and proves that when $\Lambda=\Z_p$ for a prime $p$ invertible in $\k$ we have $\mathrm{B}=\mathscr{O}(G^\vee_{\Z_p})$, which finishes the proof. 

\begin{rmk}
\phantomsection
\label{rmk:remarks-proof-Satake}
\begin{enumerate}
	\item In the course of the proof outlined above, one constructs a closed immersion of group schemes $T^\vee_\Lambda \subset \Spec(\mathrm{B})$; this construction uses the weight functors of~\S\ref{sss:weight-functors}, and more specifically the decomposition in Theorem~\ref{thm:weight-functors}\eqref{it:weight-functors-2}.
	\item The identification of $\mathrm{B}$ (in case $\Lambda=\Z_p$) with $\mathscr{O}(G^\vee_{\Z_p})$ relies on results of Prasad--Yu 
and a special analysis of the case when $\Lambda$ is an algebraically closed field; see~\cite[\S 1.14]{BR18} for details. For the application discussed below, the important case is when $\Lambda$ is a field of positive characteristic, but the only known proof of Theorem~\ref{thm:Satake} in this case requires to treat the case of $\Z_p$ in parallel. On the other hand, if one only cares about the case when $\Lambda$ is a field of characteristic $0$, the proof can be shortened quite a bit; in~\cite{BR18} the proof of this case only requires the first 9 sections. 
\end{enumerate}
\end{rmk}

\section{Representations of reductive groups and the Iwahori--Whittaker model}
\label{sec:lecture-3}

From now on, we assume that $\k$ has positive characteristic $\ell$, and we take for $\Lambda$ an algebraically closed field $\F$ of characteristic $p>0$, with $\ell \neq p$. 
In this section, we recall some standard results on representations of reductive algebraic groups over $\F$, and then we discuss the Iwahori--Whittaker model for the Satake category, following~\cite{BGMRR19}.

\subsection{Representations of $G^\vee_\F$} 


In this section we state some results from representation theory of reductive groups that will be used later. For a full treatment of these questions see~\cite[Part~II]{Jan03}. For a shorter overview, see~\cite{book}.

\subsubsection{Standard, costandard and simple representations} 
\label{sss:M-N-L}

Recall that in~\S\ref{sss:attractors-Gr} we have fixed a system of positive roots $R_+\subset R \subset \X=X^*(T)$. 
Let $R^\vee_+\subset R^\vee \subset \X^\vee = X^*(T^\vee_\F)$ be the corresponding subset of positive coroots. 
We will denote by $B^\vee_\F\subset G^\vee_\F$ the Borel subgroup containing $T^\vee_\F$ such that the $T^\vee_\F$-weights in the Lie algebra of $B^\vee_\F$ are $-R^\vee_+$. 

If $U^\vee_\F$ is the unipotent radical of $B^\vee_\F$, the composition $T^\vee_\F \hookrightarrow B^\vee_\F \twoheadrightarrow B^\vee_\F / U^\vee_\F$ is an isomorphism; hence any character of $T^\vee_\F$ extends in a unique way to a character of $B^\vee_\F$. Hence
each element $\lambda\in \X^\vee$ determines a $1$-dimensional $B^\vee_\F$-module $\F_{B^\vee_\F}(\lambda)$. We consider the associated \emph{costandard module} 
\[
\coweyl(\lambda) =\Ind^{G^\vee_\F}_{B^\vee_\F} (\F_{B^\vee_\F}(\lambda)) =\{f\in \mathscr{O}(G^\vee_\F) \mid \forall g\in G^\vee_\F,\, \forall b\in  B^\vee_\F,\ f(gb)=\lambda(b)^{-1} f(g)\}.
\]
(In~\cite{Jan03} this module is denoted $H^0(\lambda)$.)

These modules satisfy the following properties.

\begin{facts} 
\begin{enumerate}
\item 
\label{it:facts-coweyl-1}
We have $\coweyl(\lambda)\neq 0$ if and only if $\lambda\in \X^\vee_+$. 
In this case $\coweyl(\lambda)$ is finite dimensional. 
\item 
\label{it:facts-coweyl-2}
For any $\lambda\in \X^\vee_+$, there exists a unique simple submodule $\simp(\lambda)\subset \coweyl(\lambda)$. 
\item 
\label{it:facts-coweyl-3}
The assignment $\lambda\mapsto \simp(\lambda)$ induces a bijection between $\X^\vee_+$ and the set of isomorphism classes of simple objects in the category $\Rep(G^\vee_\F)$.
\end{enumerate} 
\end{facts}

For the proof of~\eqref{it:facts-coweyl-1}, see~\cite[\S II.2.1 and Proposition~II.2.6]{Jan03}. For~\eqref{it:facts-coweyl-2}, see~\cite[Corollary~II.2.3]{Jan03}. For~\eqref{it:facts-coweyl-3}, see~\cite[Proposition~II.2.4]{Jan03}.

For $\lambda\in \X^\vee_+$, we also consider the \emph{standard module} 
\[
\weyl(\lambda) :=\coweyl(-w_0\lambda)^*,
\]
where $w_0$ is the longest element in the Weyl group $W$ (for the Coxeter structure determined by our choice of Borel subgroups). 
For any $\lambda \in \X^\vee_+$ we have
\[
\simp(\lambda)^* \simeq \simp(-w_0\lambda),
\]
see~\cite[Corollary~II.2.5]{Jan03}. Hence $\weyl(\lambda)$ has a unique simple quotient, isomorphic to $\simp(\lambda)$.
We have 
\begin{equation}
\label{eqn:Hom-Weyl-coweyl}
\dim_\F \Hom(\weyl(\lambda), \coweyl(\mu))=\delta_{\lambda,\mu},
\end{equation}
see~\cite[Proposition~II.4.13]{Jan03}. Since we have a nonzero composition
\[
\weyl(\lambda) \twoheadrightarrow \simp(\lambda) \hookrightarrow \coweyl(\mu),
\]
it follows that $\simp(\lambda)$ is the image of any nonzero morphism $\weyl(\lambda) \rightarrow \coweyl(\lambda)$. 

\subsubsection{Corresponding perverse sheaves}
\label{sss:corresp-perv-sheaves}
 
 Now we come back to the geometry of the affine Grassmannian $\Gr_G$.
For any $\mu\in \X^\vee$, we have a (morphism representable by a) locally closed immersion $j_\mu\colon \Gr_{G,\mu}\rightarrow \Gr_G$, see~\S\ref{sss:schubert-cells}. 
We define the following objects in $\Perv_{\Loop^+G}(\Gr_G, \F)$:
\begin{align*}
J_*(\mu)= \pH^0 (j_{\mu*}\underline{\F}_{\Gr_{G,\mu}} [\langle 2\rho, \mu\rangle]), \\
J_!(\mu)= \pH^0 (j_{\mu!}\underline{\F}_{\Gr_{G,\mu}} [\langle 2\rho, \mu\rangle]), \\
J_{!*}(\mu)= (j_{\mu})_{!*} (\underline{\F}_{\Gr_{G,\mu}} [\langle 2\rho, \mu\rangle]).
\end{align*}
The general theory of perverse sheaves guarantees that each $J_{!*}(\mu)$ is a simple perverse sheaf, isomorphic to the image of any nonzero morphism from $J_{!}(\mu)$ to $J_{*}(\mu)$, and that moreover these objects are representatives for the isomorphism classes of simple objects in $\Perv_{\Loop^+G}(\Gr_G, \F)$. The following statement says that these objects are ``geometric incarnations'' of the representations considered in~\S\ref{sss:M-N-L}.

\begin{prop}
For any $\mu\in \X^\vee_+$ 
there are canonical isomorphisms
\[
\mathrm{Sat}(J_*(\mu))\simeq \coweyl(\mu), \quad
\mathrm{Sat}(J_!(\mu))\simeq \weyl(\mu), \quad 
\mathrm{Sat}(J_{!*}(\mu))\simeq \simp(\mu).
\]
\end{prop}

For a proof, see~\cite[Proposition~13.1]{MV07} or~\cite[\S 1.5.1]{central}.

\subsubsection{Tilting representations} 

It turns out that, in addition to~\eqref{eqn:Hom-Weyl-coweyl}, one can describe all $\Ext$-groups in $\Rep(G^\vee_\F)$ from a standard module to a costandard module:
we have 
\[
\Ext_{\Rep(G^\vee_\F)}^n( \weyl(\lambda),\coweyl(\mu))=0 \quad \text{for any $\lambda,\mu\in \X^\vee_+$ and $n>0$}.
\]
In fact, this statement\footnote{This statement is closely related to, but different from,~\cite[Proposition~II.4.13]{Jan03}, which claims a similar vanishing in the category of \emph{all} (not necessarily finite-dimensional) algebraic $G^\vee_\F$-modules.} is part of the claim that $\Rep(G^\vee_\F)$ is a \emph{highest weight category}, see~\cite[\S 2.2 in Chap.~1]{book}.
In this context, it is natural to study the \emph{tilting objects}, defined as follows. 

\begin{definition}
A representation $V\in \Rep(G^\vee_\F)$ is called \emph{tilting} if it admits both a filtration with subquotients of the form $\coweyl(\lambda)$ ($\lambda\in \X^\vee_+$) and a filtration with subquotients of the form $\weyl(\lambda)$ ($\lambda\in \X^\vee_+$). 
\end{definition}

The following properties are part of the general theory of tilting objects in highest weight categories, see~\cite[\S 5 in Appendix~A]{book}. In the context of representations of reductive groups they were first observed by Donkin; see~\cite[Chap.~II.E]{Jan03}.

\begin{facts}
\begin{enumerate}
\item If $V\in \Rep(G^\vee_\F)$ is a tilting representation, then the number $(V: \coweyl(\lambda))$ 
of occurrences of $\coweyl(\lambda)$ in a filtration with costandard subquotients does not depend on the filtration. 
In fact we have 
\[
(V: \coweyl(\lambda))=\dim_\F \Hom(\weyl(\lambda), V).
\]
\item 
If $V\in \Rep(G^\vee_\F)$ is a tilting representation, then the number $(V: \weyl(\lambda))$ 
of occurrences of $\weyl(\lambda)$ in a filtration with standard subquotients does not depend on the filtration. 
In fact we have 
\[
(V: \weyl(\lambda))=\dim_\F \Hom(V, \coweyl(\lambda)).
\]
\item 
For any $\lambda\in \X^\vee_+$, there exists a unique indecomposable tilting representation $\til(\lambda)\in \Rep(G^\vee_\F)$ such that $(\til(\lambda):\coweyl(\lambda))=1$ and $(\til(\lambda): \coweyl(\mu))=0$ unless $\mu\leq \lambda$. 
Moreover, the assignment $\lambda\mapsto \til(\lambda)$ induces a bijection 
from $\X^\vee_+$ to the set of isomorphism classes of indecomposable tilting $G^\vee_\F$-modules.
\item Any tilting representation is isomorphic to a direct sum of indecomposable tilting representations, and this decomposition is unique up to isomorphism and permutation of the factors. 
\end{enumerate}
\end{facts}

\subsubsection{Affine Weyl group} 
\label{sss:Waff}

From now on we assume that the characteristic $p$ of $\k$ is positive. 
Let
\[
W_\af=W\ltimes \Z R^\vee
\]
be the affine Weyl group. For $\lambda \in \Z R^\vee$, the corresponding element in $W_\af$ will be denoted $t_\lambda$.
We will consider the action of $W_\af$ on $V=\X^\vee\otimes_{\Z} \R$ given by 
\[
wt_\lambda\bullet v= w(v+p\lambda +\rho)-\rho, \quad 
\text{for $w\in W$, $\lambda\in \Z R^\vee$ and $v\in V$}.
\] 
It can be easily seen that this action stabilizes $\X^\vee \subset V$, and realizes $W_\af$ as a group of affine transformations of $V$ generated by the reflections $t_{n\beta^\vee} s_\beta$ for $\beta\in R$ and $n\in \Z$, where $s_\beta \in W$ is the reflection associated with $\beta$. 
The corresponding intersections of reflection hyperplanes provides a decomposition of $V$ as a disjoint union of the \emph{facets}, i.e. the nonempty subsets of the form 
\[
F=\ 
\left\{ 
	v \in V \left |\ 
    \begin{gathered} 
	\langle v+\rho, \alpha \rangle=n_\alpha p, \ \forall \alpha\in R_0(F) \\ 
	(n_\alpha-1)p< \langle v+\rho, \alpha\rangle < n_\alpha p,\ \forall \alpha\in R_1(F) 
	\end{gathered}  
\right. \right\}
\]
for some decomposition $R_+=R_0(F)\sqcup R_1(F)$ and some collection of numbers $(n_\alpha : \alpha \in R_+)$. (We insist that a subset defined in this way might be empty; we will call facet those \emph{which are nonempty}.)
A facet is called an \emph{alcove} if $R_0(F)=\varnothing$, and a \emph{wall} if $|R_0(F)|=1$. 
To a wall we can associate in a natural way a reflection $s_F\in W_\af$ which acts as the identity of this wall.

\begin{exm}
\label{ex:fund-alcove}
An important example of alcove is
the \emph{fundamental alcove}, defined as
\[
C=\{v\in V|\ 0< \langle v+\rho, \alpha\rangle <p,\ \forall \alpha\in R\}.
\]
We have $\X^\vee \cap C \neq \varnothing$ if and only if $p\geq h$, where $h$ is the Coxeter number of $G^\vee_\F$, see~\cite[Equation~(8) in~\S II.6.2]{Jan03}.
\end{exm}

The following properties are deduced from the general theory of (discrete) groups of affine transformations of an affine space, see~\cite[\S II.6.2]{Jan03}.

\begin{facts} 
\phantomsection
\label{facts:Waff}
\begin{enumerate}
	\item $\overline{C}$ is a fundamental domain for the action of $W_\af$ on $V$. 
	\item Denote by $S_\af$ the set of reflections $s_F$ where $F$ runs over the (finitely many) walls $F$ contained in $\overline{C}$. Then $(W_\af,S_\af)$ is a Coxeter system. 
\item For $v \in \overline{C}$, the stabilizer of $v$ in $W_\af$ is the subgroup generated by 
$S_{\af,v}:=\{s\in S_\af \mid \ s\bullet v= v\}$. 
\end{enumerate}
\end{facts}

In particular, from these properties we obtain a bijection 
\[
\bigsqcup_{\lambda\in \X^\vee\cap \overline{C}} W_\af/ \langle S_{\af,\lambda} \rangle \simto \X^\vee
\]
sending the class of $w$ in $W_\af/ \langle S_{\af,\lambda} \rangle$ to $w \bullet \lambda$.

Recall that given a subset $I \subset S_\af$ we can consider the (parabolic) subgroup $\langle I \rangle \subset W_\af$ generated by $I$. Then each coset in $\langle I \rangle \backslash W_\af$, resp.~in $W_\af / \langle I \rangle$, contains a unique element of minimal length (called the minimal element of that coset) and, if $\langle I \rangle$ is finite, a unique element of maximal length (called the maximal element of the coset).
For $\lambda\in \X^\vee\cap \overline{C}$, we set 
\[
W^{(\lambda)}_\af=
\left\{ 
	w\in W_\af \left |\ 
    \begin{gathered} 
	\text{$w$ is minimal in $Ww$} \\ 
	\text{$w$ is maximal in $w \langle S_{\af,\lambda}\rangle$} 
	\end{gathered}  
\right. \right\}.
\]
Then the assignment $w\mapsto w\bullet \lambda$ induces a bijection
\begin{equation}
\label{eqn:parametrization-dom-wts-block}
W^{(\lambda)}_\af\simto (W_\af\bullet \lambda)\cap \X^\vee_+,
\end{equation}
see~\cite[Chap.~1, \S 2.8]{book}. More specifically, if $w$ is maximal in $w \langle S_{\af,\lambda}\rangle$, then $w \bullet \lambda$ is dominant if and only if $w \in W^{(\lambda)}_\af$.

\begin{exm} 
If $p\geq h$ (where $h$ is as in Example~\ref{ex:fund-alcove}) we have $0\in \X^\vee \cap C$, and $W^{(0)}_\af=\{ w\in W_\af \mid w \text{ is minimal in }Ww\}$.
\end{exm}

\subsubsection{Linkage principle}
 
 The following statement was first conjectured by Verma, and then proved in increasing degrees of generality by Humphreys, Carter--Lusztig (for the group $\mathrm{GL}_n$), Jantzen, and finally Andersen. (See~\cite{RW22} for precise references.) It has had a deep influence on the modern approach to the representation theory of reductive algebraic groups.
 
\begin{thm}
\label{thm:linkage}
For $\lambda,\mu \in \X^\vee_+$, we have the implication
\[
\Ext^1_{\Rep(G^\vee_\F)}(\simp(\lambda), \simp(\mu)) \neq 0\quad \Rightarrow \quad W_\af\bullet \lambda=W_\af \bullet \mu.
\]
\end{thm} 

This theorem can be proved rather easily under mild assumptions by consideration of the central characters for the action of the Lie algebra of $G^\vee_\F$ (see~\cite[Chap.~1, \S 2.5]{book}), but the proof in full generality is much more subtle.

Theorem~\ref{thm:linkage} allows to break the category $\Rep(G^\vee_\F)$ into smaller pieces (often called \emph{blocks}, although this terminology is somehow not appropriate; see~\cite[Remark~1.1]{RW22}), as follows.
For $\lambda\in \X^\vee\cap \overline{C}$, we define $\Rep_{\lambda}(G^\vee_\F)$ as the Serre subcategory of $\Rep(G^\vee_\F)$ generated by the simple modules $\simp(\nu)$ with $\nu\in (W_\af\bullet \lambda)\cap \X^\vee_+$.  
Then it follows immediately from Theorem~\ref{thm:linkage} that we have a decomposition 
\begin{equation}
\label{eqn:decomposition-Rep-blocks}
\Rep(G^\vee_\F)= \bigoplus_{\lambda\in \X^\vee\cap \overline{C}} \Rep_{\lambda}(G^\vee_\F).
\end{equation}
(For details, see e.g.~\cite[Chap.~1, \S 2.5]{book}.)
Note that:
\begin{itemize}
	\item in view of~\eqref{eqn:parametrization-dom-wts-block}, the isomorphism classes of simple objects in $\Rep_\lambda(G^\vee_\F)$ are in a canonical bijection with $W_\af^{(\lambda)}$;
	\item for $\mu\in(W_\af\bullet \lambda)\cap \X^\vee_+$, the object $\til(\mu)$ belongs to $\Rep_\lambda(G^\vee_\F)$, which implies that for $\nu\in \X^\vee_+$ we have $(\til(\mu): \coweyl(\nu))=0$ unless $\nu\in W_\af\bullet \lambda$.
\end{itemize}

\subsubsection{Characters}
\label{sss:characters}

Recall any algebraic representation of $T^\vee_\F$ is the direct sum of its weight spaces: in other words,
for any $M\in \Rep(T^\vee_\F)$, we have 
\[
M=\bigoplus_{\mu\in \X^\vee} M_\mu \quad \text{where $M_\mu = \{v \in M \mid \forall t \in T^\vee_\F, \, t \cdot v = \lambda(t) v\}$}.
\]
We define the character of $M$ as the element 
\[
\mathrm{ch}(M)=\sum_{\mu\in \X^\vee} \dim(M_\mu)e^\mu\ \in \Z[\X^\vee].
\]
We will consider these characters in case $M$ is the restriction of a $G^\vee_\F$-module, and will write in this case $\mathrm{ch}(M)$ for $\mathrm{ch}(M_{| T^\vee_\F})$.

A crucial question in the representation theory of reductive groups is the following.
 
\begin{quest}
\label{quest:simple-char}
Compute the characters $\mathrm{ch}(\simp(\mu))$ for all $\mu\in \X^\vee_+$, or in other words compute the characters $\mathrm{ch}(\simp(w \bullet \lambda))$ for all $\lambda \in \X^\vee \cap \overline{C}$ and all $w \in W_\af^{(\lambda)}$.
\end{quest}

This question is the subject of a highly influential conjecture of Lusztig, see~\S\ref{sss:Lusztigs-formula} below or~\cite[Chap.~II.C]{Jan03}. This conjecture is now known to be true in large characteristic (i.e.~when $p$ is larger than a bound depending on the root system $R^\vee$) thanks to works of Kazhdan--Lusztig, Kashiwara--Tanisaki and Andersen--Jantzen--Soergel, but not in any satisfactory generality (due to more recent work of Williamson.) See~\cite[Chap.~1, \S 4.4]{book} for historical details. In any case, we expect an answer to this question expressing $\mathrm{ch}(\simp(w \bullet \lambda))$ in terms of some combinatorial data attached to the Coxeter system $(W_\af, S_\af)$, evaluated for $w$.

Here we are interested in a different way of formulating an answer to this question, first proposed by Andersen. Namely, thanks to results of Andersen (and a recent improvement due to Sobaje), one knows that to give an answer to this question it suffices to compute the characters 
$\mathrm{ch}(\til(\mu))$ for all $\mu\in \X^\vee_+$, see~\cite{sobaje} for details. (The concrete way to deduce simple characters from tilting characters is admitedly not easy to explain.)
Recall that the characters $\mathrm{ch}(\coweyl(\nu))$ ($\nu \in \X^\vee_+$) are known, and given by Weyl's character formula, see~\cite[Corollary~II.5.11]{Jan03}.
Therefore, the problem of computing $\mathrm{ch}(\til(\mu))$ is (in theory) equivalent to that of computing
$(\til(\mu):\coweyl(\nu))$ for all $\nu\in \X^\vee_+$. 
Hence Question~\ref{quest:simple-char} is equivalent to the following question.

\begin{quest}
\label{quest:tilting-mult}
Compute the multiplicities $(\til(\mu):\coweyl(\nu))$ for all $\mu,\nu\in \X^\vee_+$, or equivalently compute the multiplicities $\big(\til(w\bullet \lambda): \coweyl(y\bullet \lambda)\big)$ for all $\lambda\in \X^\vee \cap \overline{C}$ and all $w,y\in W_\af^{(\lambda)}$.
\end{quest}

Here again, the answer we expect should express $\big(\til(w\bullet \lambda): \coweyl(y\bullet \lambda)\big)$ in terms of some Coxeter-theoretic data attached to $w$ and $y$.

\subsection{Iwahori--Whittaker model}


In this section we explain the main result of~\cite{BGMRR19}, which provides ``another incarnation of the Satake category $\Perv_{\Loop^+ G}(\Gr_G, \F)$'', i.e.~an equivalence between $\Perv_{\Loop^+ G}(\Gr_G, \F)$ and another category of perverse sheaves on $\Gr_G$. The latter category is \emph{not} monoidal, but only a right module over $\Perv_{\Loop^+ G}(\Gr_G, \F)$; this equivalence will be an equivalence of module categories over $\Perv_{\Loop^+ G}(\Gr_G, \F)$. Nevertheless this category sometimes has some advantages over $\Perv_{\Loop^+ G}(\Gr_G, \F)$, as we will explain in~\S\ref{sss:IW-model}.

\subsubsection{Iwahori subgroup}

Consider the evaluation map $\mathrm{ev}_{z=0}\colon \Loop^+G\rightarrow G$ induced by sending $z$ to $0$, i.e.~given at the level of $R$-points by the morphism $G(R\br{z}) \to G(R)$ induced by the ring morphism $R\br{z} \to R\br{z}/(z) = R$.
The \emph{Iwahori subgroup} associated with $B^+$ is the subgroup scheme 
\[
I^+:= (\mathrm{ev}_{z=0})^{-1}(B^+)\subset \Loop^+G.
\]
We will also consider the pro-unipotent radical of $I^+$, defined as $I^+_{\un}=(\mathrm{ev}_{z=0})^{-1}(U^+)$ where, as in~\S\ref{sss:si-orbits}, $U^+$ is the unipotent radical of $B^+$. For $\lambda\in \X^\vee$, we set 
\[
\sO_\lambda=I^+ \cdot \Loop_\lambda=I^+_{\un} \cdot \Loop_\lambda,
\]
and let $j'_\lambda\colon \sO_\lambda\rightarrow \Gr$ be the embedding. Here $\sO_\lambda$ is a locally closed subscheme of the projective scheme $\Gr_{\leq \mathsf{dom}(\lambda)}$, and it is well known to be isomorphic to an affine space over $\k$.
We have 
\[
|\Gr_G|=\bigsqcup_{\lambda\in \X^\vee}|\sO_\lambda|,
\]
and, for all $\mu\in \X^\vee_+$, 
\[
|\Gr_{G,\mu}|=\bigsqcup_{\lambda\in W\mu}|\sO_\lambda|.
\]

\subsubsection{Artin--Schreier sheaf} 


Consider the morphism 
\[
\mathrm{AS}\colon \A^1_\k \rightarrow \A^1_\k
\]
given by $x\mapsto x^\ell-x$. This morphism is a Galois covering with Galois group $\F_\ell$. 
Since $\ell$ is invertible in $\F$, pushing forward the constant sheaf along this map provides a local system with a decomposition 
\[
\mathrm{AS}_*\underline{\F}_{\A^1_\k}\ =\bigoplus_{\psi\colon \F_\ell \rightarrow \F^\times} \scL_{\psi}.
\]
We fix a choice of non-trivial morphism $\psi \colon \F_\ell \rightarrow \F^\times$, and set $\scL_{\mathrm{AS}}=\scL_\psi$. Since, the local system associated with the trivial morphism $\psi$ is the constant local system, this local system must satisfy
\begin{equation}
\label{eqn:vanishing-AS}
\coH^\bullet(\A^1_\k, \scL_{\mathrm{AS}})=0=\coHc^\bullet(\A^1_\k, \scL_{\mathrm{AS}}).
\end{equation}

\begin{rmk}
In an analytic setting (as discussed in Remark~\ref{rmk:etale-const}), there exists no local system on $\A^1$ which satisfies~\eqref{eqn:vanishing-AS}. This is the main reason why we need to work with \'etale sheaves in~\cite{BGMRR19} and~\cite{RW22} (and, therefore, in these notes).
\end{rmk}

\subsubsection{Iwahori--Whittaker sheaves}
\label{sss:IW-sheaves}

For any simple root $\alpha$, we will denote by $U_\alpha \subset U^+$ the corresponding $1$-parameter subgroup. 
We consider the composition of morphisms of group schemes 
\[
I^+_{\un} \xrightarrow{\mathrm{ev}_{z=0}} U^+\rightarrow U^+/[U^+,U^+]=\prod_{\alpha\text{ simple root}} U_\alpha \rightarrow \G_{\mathrm{a},\k},
\] 
where the second morphism is the obvious quotient map, and the last morphism is a morphism which is nontrivial on each $U_\alpha$. 
Pulling back $\scL_{\mathrm{AS}}$ through this morphism we get a local system on $I^+_{\un}$, which will also be denoted $\scL_{\mathrm{AS}}$. 
Let $a\colon I^+_{\un}\times \Gr_G \rightarrow \Gr_G$ be the action map. 
We will denote by
\[
\Db_{\sI\sW}(\Gr_G,\F)
\]
the full subcategory of $\Dbc(\Gr_G,\F)$ whose objects are the complexes $\scF$ such that there exists an isomorphism 
\begin{equation}
\label{eqn:condition-AS}
a^*\scF\simto \scL_{\mathrm{AS}}\boxtimes_{\F} \scF.
\end{equation}
The same formulation as for the convolution product $\star$ in~\S\ref{sss:convolution} provides a bifunctor
\[
\Dbc(\Gr_G,\F) \times \Db_{\Loop^+ G}(\Gr_G,\F) \to \Dbc(\Gr_G,\F)
\]
which defines a right action of the monoidal category $\Db_{\Loop^+ G}(\Gr_G,\F)$ on $\Dbc(\Gr_G,\F)$. It is well known that this bifunctor is t-exact, in the sense that the image of a pair of perverse sheaves is again perverse, see~\cite[Lemma~2.3]{BGMRR19}. It is also restricts to a bifunctor
\[
\Db_{\sI\sW}(\Gr_G,\F) \times \Db_{\Loop^+ G}(\Gr_G,\F) \to \Db_{\sI\sW}(\Gr_G,\F)
\]
which defines a right action of $\Db_{\Loop^+ G}(\Gr_G,\F)$ on $\Db_{\sI\sW}(\Gr_G,\F)$.

For any $\lambda \in \X^\vee$ we can define an a similar way the full subcategory $\Db_{\sI\sW}(\sO_\lambda ,\F) \subset \Dbc(\sO_\lambda ,\F)$.
We set
\[
\X^\vee_{++}=\{ \lambda\in \X^\vee \mid \forall \alpha\in R_+,\, \langle \lambda,\alpha\rangle>0\}.
\]
The following are basic properties of these categories. For references, see~\cite[\S 5.1]{RW22}.

\begin{facts} 
\phantomsection
\label{facts:IW}
\begin{enumerate}
\item $\Db_{\sI\sW}(\Gr_G,\F)$ is a triangulated subcategory of $\Dbc(\Gr_G,\F)$. 
Moreover, the perverse t-structure on $\Dbc(\Gr_G,\F)$ restricts to a t-structure on the subcategory $\Db_{\sI\sW}(\Gr_G,\F)$. 
\item 
\label{it:IW-vanishing}
For $\lambda\in \X^\vee$ we have $\Db_{\sI\sW}(\sO_\lambda ,\F)=0$ unless $\lambda\in \X^\vee_{++}$. 
\item Let $\lambda\in \X^\vee_{++}$. 
There exists a unique (up to isomorphism) local system $\scL_{\mathrm{AS},\lambda}$ on $\sO_\lambda$ whose pullback under the morphism $I_{\un}^+\rightarrow \sO_\lambda$ given by $g\mapsto g \cdot \Loop_\lambda$ is $\scL_{\mathrm{AS}}$. 
\end{enumerate}
\end{facts} 

The proof of~\eqref{it:IW-vanishing} is based on the observation that if $\lambda \in \X^\vee \smallsetminus \X^\vee_{++}$, the stabilizer of $\Loop_\lambda$ in $I^+_{\un}$ contains a simple root subgroup $U_\alpha$, so that the orbit $\sO_\lambda$ cannot support any nonzero complex $\scF$ which satisfies~\eqref{eqn:condition-AS}. This vanishing implies that all complexes in $\Db_{\sI\sW}(\Gr_G,\F)$ have trivial restrictions and corestrictions to all $\sO_\lambda$ with $\lambda \in \X^\vee \smallsetminus \X^\vee_{++}$. In other words, the category $\Db_{\sI\sW}(\Gr_G,\F)$ ``does not see these orbits.''

For any $\lambda\in \X^\vee_{++}$, we set 
\begin{align*}
\Delta^{\sI\sW}_\lambda &=(j'_\lambda)_! \scL_{\mathrm{AS},\lambda}[\dim \sO_\lambda], \\
\nabla^{\sI\sW}_\lambda &=(j'_\lambda)_* \scL_{\mathrm{AS},\lambda}[\dim \sO_\lambda], \\
\IC^{\sI\sW}_\lambda &=(j'_\lambda)_{!*} (\scL_{\mathrm{AS},\lambda}[\dim \sO_\lambda]).
\end{align*}
Here, in constrast with the situation considered in~\S\ref{sss:corresp-perv-sheaves}, we do not need to take perverse cohomology in the definition of $\Delta^{\sI\sW}_\lambda$ and $\nabla^{\sI\sW}_\lambda$: these complexes are perverse sheaves because $j'_\lambda$ is an affine morphism. These objects can be used to characterize the restriction of the perverse t-structure to $\Db_{\sI\sW}(\Gr_G,\F)$:
we have 
\begin{align*}
{}^{\mathrm{p}} \hspace{-1pt} D^{\leq 0}_{\sI\sW}(\Gr_G,\F) & =\langle \Delta^{\sI\sW}_\lambda[n] : \lambda\in \X^\vee_{++},\, n\in \Z_{\geq 0}\rangle_{\mathrm{ext}}, \\
{}^{\mathrm{p}} \hspace{-1pt} D^{\geq 0}_{\sI\sW}(\Gr_G,\F) &= \langle \nabla^{\sI\sW}_\lambda[n] : \lambda\in \X^\vee_{++},\ n\in \Z_{\leq 0}\rangle_{\mathrm{ext}}
\end{align*}
where $\langle - \rangle_{\mathrm{ext}}$ means the full subcategory generated under extensions (in the sense of triangulated categories) by the collection of objects under consideration.
As in the setting of~\S\ref{sss:corresp-perv-sheaves},
the assignment $\lambda\mapsto \IC^{\sI\sW}_\lambda$ induces a bijection between $\X^\vee_{++}$ and the set of isomorphism classes of simple objects in $\Perv_{\sI\sW}(\Gr_G,\F)$, and there exist natural morphisms 
\[
\Delta^{\sI\sW}_\lambda\twoheadrightarrow \IC^{\sI\sW}_\lambda \hookrightarrow \nabla^{\sI\sW}_\lambda.
\]

\subsubsection{Iwahori--Whittaker model}
\label{sss:IW-model}

From now on we will assume that there exists $\xi\in \X^\vee$ such that $\langle \xi, \alpha \rangle =1$ for all simple root $\alpha$. This assumption holds at least if the quotient $\X / \Z R$ is torsion-free, e.g.~if $G$ is semisimple of adjoint type; this is a very mild assumption since essentially any question one can ask about the representation theory of $G^\vee_\F$ can be reduced to the case this condition is satisfied.
We fix such a $\xi$; then $\sO_\xi$ is minimal (for the order given by inclusions of closures) among all the orbits $\{\sO_\lambda\}_{\lambda\in \X^\vee_{++}}$, so that we have isomorphisms 
\[
\Delta^{\sI\sW}_\xi=\IC^{\sI\sW}_\xi=\nabla^{\sI\sW}_\xi.
\]

The following statement is the main result of~\cite{BGMRR19}.

\begin{thm}
\label{thm:IW-model}
The functor $\scF \mapsto \Delta^{\sI\sW}_\xi \star \scF$ 
induces an equivalence of categories
\[
\Perv_{\Loop^+G}(\Gr_G,\F) \simto \Perv_{\sI\sW}(\Gr_G,\F)
\]
which, for any $\lambda\in \X^\vee_+$,
sends the object $J_!(\lambda)$, resp.~$J_*(\lambda)$, resp.~$J_{!*}(\lambda)$, to $\Delta^{\sI\sW}_{\lambda+\xi}$, resp.~$\nabla^{\sI\sW}_{\lambda+\xi}$, resp.~$\IC^{\sI\sW}_{\lambda+\xi}$.
\end{thm}

\begin{rmk}
A version of this statement for
characteristic-$0$ coefficients was proved earlier by Arkhipov--Bezru\-kavnikov--Braverman--Gaitsgory--Mirkovi\'{c}, see~\cite{ABBGM05}. 
In this case, both categories in this theorem are semisimple. 
\end{rmk} 

Theorem~\ref{thm:IW-model} should be thought of as giving a different incarnation of the Satake category $\Perv_{\Loop^+G}(\Gr_G,\F)$, where the tensor structure is not visible but which has other avantages. In particular, being equivalent to $\Rep(G^\vee_\F)$, the category $\Perv_{\Loop^+G}(\Gr_G,\F)$ is a highest weight category. This property is in fact (implicitly) proved by Mirkovi{\'c}--Vilonen in the course of their proof of Theorem~\ref{thm:Satake}, see~\cite[Proposition~1.12.4]{BR18}, but this proof relies on a very clever use of some specific properties of the geometry of $\Gr_G$. On the other hand, the fact that the category $\Perv_{\sI\sW}(\Gr_G,\F)$ has a highest weight structure has a transparent explanation: using standard results of Be{\u\i}linson--Ginzburg--Soergel, this is simply a consequence of the fact that the $I^+_{\un}$-orbits on $\Gr_G$ are affine spaces, see~\cite[Corollary~3.6]{BGMRR19}. This observation is crucial for the applications discussed below.

\subsubsection{Application to tilting representations} 

Let us now consider the category
\[
\Db_{(\Loop^+G)}(\Gr_G,\F)
\]
as the category of ``$\Loop^+G$-constructible'' (rather than equivariant) complexes on $\Gr_G$, i.e.~the full triangulated subcategory of $\Dbc(\Gr_G,\F)$ generated by the objects $J_{!*}(\lambda)$ for $\lambda\in \X^\vee_+$. 
The perverse t-structure restricts to a t-structure on $\Db_{(\Loop^+G)}(\Gr_G,\F)$, whose heart is denoted by $\Perv_{(\Loop^+G)}(\Gr_G,\F)$. We also have a natural t-exact ``forgetful'' functor
\begin{equation}
\label{eqn:For-DbGr}
\Db_{\Loop^+G}(\Gr_G,\F) \to \Db_{(\Loop^+G)}(\Gr_G,\F).
\end{equation}

The following result is due to Mirkovi\'{c}--Vilonen; see~\cite[Proposition~2.1]{MV07} or~\cite[\S 1.10.2]{BR18}.

\begin{prop}
\label{prop:For-equiv}
The functor~\eqref{eqn:For-DbGr} restricts to an equivalence of categories 
\[
\Perv_{\Loop^+G}(\Gr_G,\F)\simto \Perv_{(\Loop^+G)}(\Gr_G,\F).
\]
\end{prop}

\begin{rmk}
Let us insist that the functor~\eqref{eqn:For-DbGr} is far from being fully faithful; only its restriction to perverse sheaves is. The fully faithfulness statement in Proposition~\ref{prop:For-equiv} follows from the general theory of perverse sheaves, but the fact that this functor is essentially surjective uses some specific geometric features of $\Gr_G$.
\end{rmk}

The following definition is due to Juteau--Mautner--Williamson~\cite{JMW14}. (This definition makes sense, and is useful, in a larger geometric setting; here we specialize it to the setting we will require.)

\begin{definition}
\label{def:parity}
A complex $\scF\in \Db_{(\Loop^+G)}(\Gr_G,\F)$ is called even (resp. odd) if
\[
\mathscr{H}^n(j^{*}_\lambda \scF)=0=\mathscr{H}^n(j^{!}_\lambda \scF)
\]
for all $\lambda\in \X^\vee_+$ unless $n$ is even (resp. odd). 
A \emph{parity complex} is a direct sum of an even and an odd object. 
\end{definition}

\begin{rmk}
The definition of an even complex can equivalently be stated as requiring that $\mathscr{H}^n(\scF)=0=\mathscr{H}^n(\mathbb{D}(\scF))$ unless $n$ is even, where $\mathbb{D}$ is the Verdier duality functor.
\end{rmk}

The following theorem is the special case, in our current setting, of the classification theorem for parity complexes established by Juteau--Mautner--Williamson in~\cite{JMW14}.

\begin{thm}
\label{thm:classification-parity-1}
For any $\lambda\in \X^\vee_+$, there exists a unique indecomposable parity complex $\scE_\lambda\in \Db_{(\Loop^+G)}(\Gr_G,\F)$ supported on $\Gr_{G,\leq \lambda}$ and such that $j^*_\lambda \scE_\lambda = \underline{\F}_{\Gr_{G,\lambda}}[\langle 2\rho, \lambda\rangle]$. 
Moreover, the assignment $(\lambda,n) \mapsto \scE_\lambda[n]$ induces a bijection between $\X^\vee_+ \times \Z$ and the set of isomorphism classes of indecomposable parity complexes in $\Db_{(\Loop^+G)}(\Gr_G,\F)$. 
\end{thm}

\begin{rmk}
\phantomsection
\label{rmk:parities-Gr}
\begin{enumerate}
\item
In the general geometric setting studied in~\cite{JMW14}, there is a \emph{unicity} statement for parity complexes, see~\cite[Theorem~2.12]{JMW14}, but the object associated with some strata might not exist, see~\cite[\S 2.3.4]{JMW14}. This problem does not occur for Kac--Moody (possibly partial) flag varieties or (possibly parahoric) affine flag varieties, by the considerations in~\cite[\S 4.1]{JMW14}.
\item
\label{it:parities-Gr-char0}
Definition~\ref{def:parity} also makes sense for coefficients $\Q_p$ rather than $\F$, and Theorem~\ref{thm:classification-parity-1} also holds in this setting. In this case, the analogue of $\scE_\lambda$ turns out to be the intersection cohomology complex associated with the constant local system on the orgit $\Gr_{G,\lambda}$. This description no longer holds for positive-characteristic coefficients.
\end{enumerate}
\end{rmk}

The following statement makes a connection between parity complexes and tilting representations.

\begin{thm}
\phantomsection
\label{thm:parity-tilting}
\begin{enumerate}
\item 
\label{it:parity-tilting-1}
The inverse images under $\mathrm{Sat}$ of the tilting $G^\vee_\F$-representations are the direct sums of direct summands of objects $\pH^0(\scE_\lambda)$ for $\lambda\in \X^\vee_+$. 
\item 
\label{it:parity-tilting-2}
If $p$ is good for $G$, then each $\scE_\lambda$ is perverse, and we have $\mathrm{Sat}(\scE_\lambda)=\til(\lambda)$ for any $\lambda\in \X^\vee_+$. 
\end{enumerate}
\end{thm}

Both statements in this theorem were conjectured by Juteau--Mautner--Williamson, and proved by them under some technical assumptions on $p$ in~\cite{JMW16}. The general case of~\eqref{it:parity-tilting-1} was proved in~\cite[\S 4.3]{BGMRR19}, as an application of Theorem~\ref{thm:IW-model}. The general case of~\eqref{it:parity-tilting-2} is proved in~\cite{mautner-riche}. (It is known that this statement does not hold in general if $p$ is bad; see~\cite[\S 3.5]{JMW16}.)

\begin{rmk}
\begin{enumerate}
	\item 
	In Theorem~\ref{thm:parity-tilting} we implicitly use Theorem~\ref{prop:For-equiv} to ``lift'' perverse sheaves in $\Db_{(\Loop^+G)}(\Gr_G,\F)$ to objects of $\Perv_{\Loop^+G}(\Gr_G,\F)$.
	\item The prime numbers which are bad (i.e.~not good) are given by the following table in terms of simple constituents in the root system of $G$ (or of $G^\vee_\F$):
\[
\begin{array}{|c|c|c|c|} 
\hline
\mathsf{A}_n & \mathsf{B}_n, \mathsf{C}_n, \mathsf{D}_n & \mathsf{E}_6, \mathsf{E}_7, \mathsf{F}_4, \mathsf{G}_2 & \mathsf{E}_8 \\ \hline 
\varnothing & 2 & 2,3 & 2,3,5  \\
\hline
\end{array}
\]
\item
From Theorem~\ref{thm:parity-tilting}\eqref{it:parity-tilting-1} one can obtain a simple proof of an important result due to Donkin in most cases and Mathieu in full generality, namely that the tensor product of two tilting $G^\vee_\F$-modules is again tilting. For details, see~\cite[Theorem~4.16]{BGMRR19}.
\end{enumerate}
\end{rmk}

\section{Smith--Treumann theory and tilting character formula} 
\label{sec:lecture-4}

In this section, we explain how to use Smith--Treumann theory on the affine Grassmannian and the geometric Satake equivalence (Theorem~\ref{thm:Satake}) to provide some answer to Question~\ref{quest:tilting-mult}.
The main reference for this section is~\cite{RW22}. 

\subsection{Preliminaries} 

\subsubsection{Smith--Treuman theory} 

\emph{Smith--Treuman theory} originates in the work of Smith in topology from 1930's, and a re-interpretation of these constructions in a more sheaf-theoretic context by D.~Treumann;\footnote{Treumann works in an analytic setting, while we rather want to work with \'etale sheaves. But the two versions have similar properties.} see~\cite{RW22} for references.

The geometric setting is the following. Let $X$ be a separated $\k$-scheme of finite type endowed with an action\footnote{For technical reasons we will assume that this action is admissible in the sense recalled in~\cite[\S 2.3]{RW22}. This condition is automatic if $X$ is quasi-projective over $\k$.} of $\G_{\mathrm{m},\k}$. 
We consider as in Section~\ref{sec:lecture-3} \'etale sheaves with coefficients in $\F$, an algebraically closed field of characteristic $p$ which is invertible in $\k$.
Let $\mu_p\subset \G_{\mathrm{m},\k}$ be the subgroup of $p$-th root of unity (a smooth constant group scheme over $\k$), and consider the embedding $i\colon X^{\mu_p}\hookrightarrow X$. (Here, $X^{\mu_p}$ denotes the closed subscheme of fixed points for the action.)

Then if $\dag \in \{*,!\}$
we can consider the following composition, where the right-hand equivalence follows from standard results on equivariant derived categories: 
\[
\Db_{\G_{\mathrm{m},\k}}(X,\F)\xrightarrow{\Res^{\G_{\mathrm{m},\k}}_{\mu_p}} \Db_{\mu_p}(X,\F) 
\xrightarrow{i^\dag} \Db_{\mu_p}(X^{\mu_p},\F) 
\simeq \Dbc(X^{\mu_p},\F[\mu_p]).
\]
(Here, $\Dbc(X^{\mu_p},\F[\mu_p])$ is the bounded derived category of sheaves of $\F[\mu_p]$-modules on $X^{\mu_p}$ with constructible cohomology.)

Recall that a complex of modules over a ring is called \emph{perfect} if it is isomorphic, in the derived category of modules over that ring, to a bounded complex of projective modules. In case the ring has finite homological dimension this definition is pointless (it is equivalent to requiring that the complex be bounded), but when this is not the case this is an important notion. One of the simplest examples of a ring which is \emph{not} of finite homological dimension is $\F[\mu_p]$; for this example, it is a standard observation that the trivial module $\F$ (i.e.~the module associated with the trivial representation of $\mu_p$), seen as a complex concentrated in degree $0$, is not perfect.

\begin{definition}
A complex $\scF\in \Db_{\mu_p}(X^{\mu_p},\F)=\Dbc(X^{\mu_p},\F[\mu_p])$ is said to \emph{have perfect geometric stalks}, if for any geometric point $\overline{x}$ of $X^{\mu_p}$ the complex $\scF_{\overline{x}}\in \Db(\F[\mu_p]\mhyphen\mathrm{Mod})$ is perfect.
\end{definition}

We define the \emph{Smith category}
\begin{equation}
\label{eqn:Smith-cat}
\mathrm{Sm}(X^{\mu_p},\F)
\end{equation}
 as the (Verdier) quotient of $\Db_{\G_{\mathrm{m},\k}}(X^{\mu_p},\F)$ by the full subcategory of complexes whose images in $\Db_{\mu_p}(X^{\mu_p},\F)$ have perfect geometric stalks. This category is triangulated, but it really looks different from derived categories of abelian categories; in fact, in this category we have $[2]=\id$. 
 
 \begin{rmk}
 A notable difference between our considerations on Smith--Treumann theory and earlier treatments in the literature is that we define $\mathrm{Sm}(X^{\mu_p},\F)$ as a quotient of the $\G_{\mathrm{m},\k}$-equivariant derived category rather than the $\mu_p$-equivariant version. This is required for our considerations of parity complexes in this category.
 \end{rmk}

Here comes a crucial observation of Treumann. 

\begin{prop} 
\label{prop:smith-restriction}
The compositions 
\[
\Db_{\G_{\mathrm{m},\k}}(X,\F)
\xrightarrow{i^*}
\Db_{\G_{\mathrm{m},\k}}(X^{\mu_p},\F) \rightarrow \mathrm{Sm}(X^{\mu_p},\F)
\]
and
\[
\Db_{\G_{\mathrm{m},\k}}(X,\F)
\xrightarrow{i^!}
\Db_{\G_{\mathrm{m},\k}}(X^{\mu_p},\F) \rightarrow \mathrm{Sm}(X^{\mu_p},\F)
\]
(where, in each case, the second functor is the natural quotient functor)
are canonically isomorphic to each other. 
\end{prop}

The functor of Proposition~\ref{prop:smith-restriction} will be denoted $i^{!*}$. 
This functor ``commutes with all sheaf operations'' in a sense similar to that discussed in Remark~\ref{rmk:HL}\eqref{it:HL-commutes}, because it can be written both a $*$- and as a $!$-pullback. 


The rough idea of the proof of Proposition~\ref{prop:smith-restriction} is to express the difference between the two functors in terms of a complex on $X \smallsetminus X^{\mu_p}$, where the action is free, which forces this difference to involve only perfect complexes.

\subsubsection{Bruhat--Tits theory}

We now consider the reductive group
\[
G\times_{\Spec(\k)} \Spec(\k\pa{z^p})
\]
over $\k\pa{z^p}$ with its maximal torus $T\times_{\Spec(\k)} \Spec(\k\pa{z^p})$. 
The apartment in the associated Bruhat--Tits building identifies with $V=\X^\vee\otimes_\Z \R$. 
We have a hyperplane arrangement in $V$, hence a system of facets, which is the same as that considered in~\S\ref{sss:Waff}, up to a shift by $\rho$. 
In particular, we have an alcove 
\[
\mathbf{a}_0=\{v \in V \mid \forall \alpha\in R_+,\, -p< \langle v,\alpha\rangle <0\}.
\]
For each facet $\mathbf{f} \subset \overline{\mathbf{a}_0}$, Bruhat--Tits theory provides us with a parahoric subgroup scheme $P_{\mathbf{f}}$ over $\Spec (\k\br{z^p})$ with generic fiber $G\times_{\Spec(\k)} \Spec(\k\br{z^p})$. 
We can consider the associated affine Grassmannian 
\[
\GGr_{P_{\mathbf{f}}}= \bigl( \Loop_p G/\L_p^+ P_{\mathbf{f}} \bigr)_{\mathrm{\acute{e}t}},
\]
where $\Loop_p G$ is loop group with $z$ replaced by $z^p$, and similarly $\L_p^+$ is defined as for $\L^+$, but for group schemes over $\k\br{z^p}$.

\begin{exm} 
\phantomsection
\label{ex:Gr-Bruhat-Tits}
\begin{enumerate}
\item
\label{it:FlGp}
If $\mathbf{f}=\mathbf{a}_0$, then $P_{\mathbf{a}_0}$ is the Iwahori group scheme associated with $B^-$, and we have 
\[
\L_p^+ P_{\mathbf{a}_0}=I^-_p=(\mathrm{ev}_{z^p=0})^{-1}(B^-)\subset \Loop^+_p G.
\] 
(Here, $\mathrm{ev}_{z^p=0}$ is similar to $\mathrm{ev}_{z=0}$, with $z$ replaced by $z^p$. The corresponding affine Grassmannian will be denoted $\Fl_{G,p}$ in this case. (This ind-scheme is often called the \emph{affine flag variety}.)
\item
If $\mathbf{f}=\{\nu\in V \mid \forall \alpha\in R,\, \langle\nu,\alpha\rangle=0\}$, then we have $P_{\mathbf{f}} = G\times_{\Spec(\k)} \Spec(\k\br{z^p})$, hence $\L_p^+ P_{\mathbf{f}}=\Loop^+_p G$. The corresponding affine Grassmannian will be denoted $\Gr_{G,p}$ in this case; it is isomorphic to $\Gr_G$.
\end{enumerate}
\end{exm}

In general if $\mathbf{f}'\subset \overline{\mathbf{f}}$, then we have a natural closed immersion $P_{\mathbf{f}}\subset P_{\mathbf{f}'}$, hence a canonical morphism $\GGr_{P_{\mathbf{f}}}\rightarrow \GGr_{P_{\mathbf{f}'}}$. 

\begin{rmk}
From Bruhat--Tits theory we have an action of $W_\af$ on $V$. 
If we define an action $\square_p$ by setting, for $w\in W$, $\mu\in \Z R^\vee$, and $v \in V$,
\[
(wt_\mu) \square_p v = w(v+p\mu),
\]
then the action of $w$ is given by $v \mapsto -w\square_p(-v)$. 
\end{rmk}

\subsubsection{$p$-canonical basis} 

Recall (from Facts~\ref{facts:Waff}) that the pair $(W_\af, S_\af)$ is a Coxeter system. As such, there is an associated Hecke algebra, called the \emph{affine Hecke algebra}, and denoted $\sH_\af$. This is a $\Z[v^{\pm 1}]$-algebra, which is free with basis $\{H_w :  w\in W_\af\}$ over $\Z[v^{\pm 1}]$, and with multiplication determined by the following rules:
\begin{align*}
&(H_s+v)(H_s-v^{-1})=0 \quad \text{for all $s\in S_\af$}, \\ 
&H_w \cdot H_y=H_{wy}, \quad \text{for all $w,y \in W_\af$ such that $\ell(wy)=\ell(w)+\ell(y)$}.
\end{align*} 

Consider the \emph{affine flag variety} 
$\Fl_{G,p}$ from Example~\ref{ex:Gr-Bruhat-Tits}\eqref{it:FlGp},
and denote by $\Fl_{G,p}^0$ the connected component of the base point.
There is a canonical action of $I^-_p$ on $\Fl_{G,p}$, induced by multiplication on the left on $\Loop_p G$; this action stabilizes $\Fl_{G,p}^0$, and the orbits on this sub-ind-scheme are parametrized in a natural way by $W_\af$: we will denote by
\[
\Fl_{G,p}^0 = \bigsqcup_{w\in W_\af} \Fl_{G,p,w}
\]
the corresponding decomposition into orbits. (Here, each $\Fl_{p,w}$ is a scheme, isomorphic to an affine space of dimension $\ell(w)$.)
In the $I^-_p$-equivariant derived category $\Db_{I^-_p}(\Fl^0_{G,p},\F)$ of sheaves on $\Fl_{G,p}$, we have a theory of parity complexes, defined by an obvious variant of Definition~\ref{def:parity}, and a classification of indecomposable parity complexes similar to Theorem~\ref{thm:classification-parity-1}, as follows.

\begin{thm}
\label{thm:classification-parity-2}
For any $w \in W_\af$, there exists a unique indecomposable parity complex $\scE_w \in \Db_{I_p^-}(\Fl^0_{G,p},\F)$ supported on $\overline{\Fl_{G,p,w}}$ and whose restriction to $\Fl_{G,p,w}$ is $\underline{\F}_{\Fl_{G,p,w}}[\ell(w)]$. 
Moreover, the assignment $(w,n) \mapsto \scE_w[n]$ induces a bijection between $W_\af \times \Z$ and the set of isomorphism classes of indecomposable parity complexes in $\Db_{I_p^-}(\Fl^0_{G,p},\F)$. 
\end{thm}

For $w \in W_\af$ we set 
\begin{equation}
\label{eqn:formula-pHw}
{}^p \hspace{-1pt} \underline{H}_w=\sum_{\substack{y\in W_\af,\\ n\in \Z}} \rk \, \mathscr{H}^{-l(y)-n}(\scE_w{}_{|\Fl_{G,p,y}})
\cdot v^n\cdot H_y.
\end{equation}
Then $( {}^p \hspace{-1pt} \underline{H}_w : w\in W_\af )$ is a $\Z[v^{\pm 1}]$-basis of $\sH_\af$, called the \emph{$p$-canonical basis}. 
The \emph{$p$-Kazhdan--Lusztig polynomials} $( {}^p \hspace{-1pt} h_{y,w} : y,w\in W_\af )$ are defined by the following formula:
\[
{}^p \hspace{-1pt} \underline{H}_w=\sum_{y\in W_\af} {}^p \hspace{-1pt} h_{y,w}\cdot H_y.
\]

\begin{rmk}
As in Remark~\ref{rmk:parities-Gr}\eqref{it:parities-Gr-char0}, Theorem~\ref{thm:classification-parity-2} also holds for coefficients $\Q_p$. In this case, the counterpart of $\scE_w$ is the intersection cohomology complex associated with the constant local system on $\Fl_{G,p,w}$, and the formula in~\eqref{eqn:formula-pHw} provides the element $\underline{H}_w$ of the Kazhdan--Lusztig basis of $\sH_\af$ associated with $w$. For coefficients $\F$ this is no longer true, but the following properties holds.
\begin{enumerate}
\item
For any $w \in W_\af$, the coefficients of the expansion of ${}^p \hspace{-1pt} \underline{H}_w$ in the Kazhdan--Lusztig basis $(\underline{H}_w : w \in W_\af)$ belong to $\Z_{\geq 0}[v^{\pm 1}]$.
\item
For any fixed $w \in W_\af$, there exists a bound $N(w)$ such that for any prime $p \geq N(w)$ we have ${}^p \hspace{-1pt} \underline{H}_w = \underline{H}_w$.
\item
There exists an algorithm, implemented by Gibson, Jensen and Williamson, that performs the computation of the elements ${}^p \hspace{-1pt} \underline{H}_w$ (in the more general setting of crystallographic Coxeter groups); see~\cite{gjw}.
\end{enumerate}
It is very difficult (but, also, very important) to say anything more specific about this basis. For a more systematic study of this subject, see~\cite[Chap.~2, \S 2.13--2.14]{book}.
\end{rmk}

If $w,y\in W_\af$ are minimal in $Ww, Wy$ respectively, then we also set 
\[
{}^p \hspace{-1pt} n_{y,w}=\sum_{x\in W}(-1)^{\ell(x)}\cdot {}^p \hspace{-1pt} h_{xy,w}.
\]
These polynomials are called the \emph{antispherical $p$-Kazhdan--Lusztig polynomials} attached to $W_\af$ and $W$. (They can be interpreted as coefficients of the expansion of a $p$-canonical basis in a standard basis in the antispherical $\sH_\af$-module.)

\subsection{Smith--Treumann theory on the affine Grassmannian} 

\subsubsection{Fixed points}

There exists a natural $\G_{\mathrm{m},\k}$-action on $\Loop G$ by ``loop rotation'', i.e.~res\-caling of $z$. 
This action stabilizes $\Loop^+G$, hence provides an action on $\Gr_G$. 
We consider the action of $\mu_p\subset \G_{\mathrm{m},\k}$; for this action we have $(\Loop G)^{\mu_p}=\Loop_p G$. 

The following result is due to Bezrukavnikov, and is crucial for what follows. For a proof, see~\cite[Proposition~4.5]{RW22}.

\begin{thm}
\label{thm:fixed-pts-Gr}
We have 
\[
(\Gr_G)^{\mu_p}=\bigsqcup_{\lambda\in -\overline{\mathbf{a}_0}\cap \X^\vee} \Gr^{(\lambda)}_G,
\quad \text{with} \quad \Gr^{(\lambda)}_G=\Loop_pG\cdot \Loop_\lambda\simeq \GGr_{P_{\mathbf{f}_\lambda}}^\circ,
\]
where $\GGr_{P_{\mathbf{f}_\lambda}}^\circ$ is the connected component of the base point in $\GGr_{P_{\mathbf{f}_\lambda}}$, and $\mathbf{f}_\lambda\subset \overline{\mathbf{a}_0}$ is the facet containing $-\lambda$. 
\end{thm}

\begin{exm}
Consider the case $G=\mathrm{PGL}_2$. Then there are canonical identifications $\X^\vee=\Z$, $V=\R$, and we have $\mathbf{a}_0=(-p,0)$. 
Hence $(\Gr_{G})^{\mu_p}$ has $p+1$ connected components, $p-1$ of which are isomorphic to $\Fl_{G,p}$ and $2$ of which are isomorphic to $\Gr_p$. 
More precisely, for $\lambda=0$ the component $\Gr_G^{(0)}$ identifies with $\Loop_pG/\Loop^+_pG$, and for $\lambda=p$  the component $\Gr_G^{(p)}$ identifies with
\[
\Loop_p G \left/ 
\begin{pmatrix}
z^p & 0\\ 0 & 1
\end{pmatrix} \cdot
\Loop^+_p G \cdot
\begin{pmatrix}
z^{-p} & 0\\ 0 & 1
\end{pmatrix}. \right.
\]
\end{exm}

Consider the group scheme $I^+$, and its action on $\Gr_G$. We have
\[
(I^+)^{\mu_p} = (\mathrm{ev}_{z^p=0})^{-1}(B^+);
\]
this group scheme will be denoted $I^+_p$. 
For any $\nu\in \X^\vee$, we have $(\sO_\nu)^{\mu_p}=I^+_p\cdot \Loop_\nu$, see~\cite[Lemma~4.4]{RW22}.
Hence we have a decomposition
\[
(\Gr_G)^{\mu_p}=\bigsqcup_{\nu\in \X^\vee} (\sO_\nu)^{\mu_p}= \bigsqcup_{\nu\in \X^\vee} I^+_p\cdot \Loop_\nu.
\] 
On the other hand, the $I^+_p$-orbits on $\GGr_{P_{\mathbf{f}_\lambda}}^\circ$ are naturally parametrized by the quotient $W_\af/W_\af^{\mathbf{f}_\lambda}$, where $W_\af^{\mathbf{f}_\lambda} \subset W_\af$ is the parabolic subgroup generated by the simple reflections $s$ satisfying $s\square_p \lambda=\lambda$: this decomposition will be written as
\[
\GGr_{P_{\mathbf{f}_\lambda}}^\circ
=\bigsqcup_{w\in W_\af/W_\af^{\mathbf{f}_\lambda}} \GGr_{P_{\mathbf{f}_\lambda,w}}^\circ.
\]
In the decomposition of Theorem~\ref{thm:fixed-pts-Gr}, these orbits are related by 
\[
\GGr_{P_{\mathbf{f}_\lambda,w}}^\circ \leftrightarrow I^+_p\cdot \Loop_{w\square_p \lambda}.
\]
Below we will consider ``Iwahori--Whittaker sheaves'' with respect to the action of $I^+_p$, i.e.~complexes which satisfy the analogue of~\eqref{eqn:condition-AS} where $I^+$ is replaced by $I_p^+$.
The orbit $I^+_p\cdot \Loop_{w\square_p \lambda}$ supports a nonzero Iwahori--Whittaker local system if and only if $w\square_p \lambda\in \X^\vee_{++}$. If $w$ is assumed to be maximal in the coset $w W_\af^{\mathbf{f}_\lambda}$, this is equivalent to the condition that $w\in W^{(\lambda-\xi)}_\af$ (see~\eqref{eqn:parametrization-dom-wts-block}).

\subsubsection{Application to the Iwahori--Whittaker model}

Recall the Iwahori--Whittaker derived category $\Db_{\sI\sW}(\Gr_G,\F)$ from~\S\ref{sss:IW-sheaves}. 
The subcategory $\Perv_{\sI\sW}(\Gr_G,\F)$ has a highest weight structure with standard, resp.~costandard, objects the perverse sheaves $\Delta^{\sI\sW}_\lambda$, resp.~$\nabla^{\sI\sW}_\lambda$ ($\lambda\in \X^\vee_{++}$), see~\S\ref{sss:IW-model}.
So we can consider the tilting objects in this category, i.e.~the objects which admits both a filtration with standard subquotients, and a filtration with costandard subquotients. The general theory of highest weight categories provides a parametrization of the isomorphism classes of tilting objects by $\X^\vee_{++}$: the object associated with $\lambda$ will be denoted $\scT^{\sI\sW}_\lambda$.
Under the geometric Satake equivalence (Theorem~\ref{thm:Satake}) and the Iwahori--Whittaker model (Theorem~\ref{thm:IW-model})
\[
\Rep(G^\vee_\F)\simto \Perv_{\Loop^+G}(\Gr_G,\F)\simto \Perv_{\sI\sW}(\Gr_G,\F),
\]
we have the folowing identifications for $\lambda\in \X^\vee_+$:
\begin{gather*}
\weyl(\lambda)\mapsto J_!(\lambda) \mapsto \Delta^{\sI\sW}_{\lambda+\xi}, \\
\coweyl(\lambda)\mapsto J_*(\lambda) \mapsto \nabla^{\sI\sW}_{\lambda+\xi}.
\end{gather*}
Hence the object $\til(\lambda) \in \Rep(G^\vee_\F)$ corresponds to $\scT^{\sI\sW}_{\lambda+\xi} \in \Perv_{\sI\sW}(\Gr_G,\F)$. 

Consider the category 
\[
\Db_{\sI\sW,\G_{\mathrm{m},\k}}(\Gr_G,\F)
\]
consisting of objects in the $\G_{\mathrm{m},\k}$-equivariant derived category of sheaves on $\Gr_G$ which satisfy~\eqref{eqn:condition-AS}. As in Facts~\ref{facts:IW}, the category $\Db_{\sI\sW,\G_{\mathrm{m},\k}}(\Gr_G,\F)$ has a natural triangulated structure, and a natural perverse t-structure. The heart of this t-structure will be denoted $\Perv_{\sI\sW,\G_{\mathrm{m},\k}}(\Gr_G,\F)$.

The following lemma is not difficult to check, see~\cite[Theorem~5.2]{RW22}.

\begin{lem}
The forgetful functor 
\[
\Perv_{\sI\sW,\G_{\mathrm{m},\k}}(\Gr_G,\F)\rightarrow \Perv_{\sI\sW}(\Gr_G,\F)
\]
is an equivalence of categories. 
\end{lem}

We set $I^+_{p,\un}=I^+_p\cap I^+_{\un}$. 
Then we have natural morphisms $I^+_{p,\un} \hookrightarrow I^+_{\un} \rightarrow U^+$. We can play the same game as in the definition of $\Db_{\sI\sW,\G_{\mathrm{m},\k}}(\Gr_G,\F)$, now the action of $I^+_{p,\un}$ on $(\Gr_G)^{\mu_p}$; the corresponding category will be denoted
\[
\Db_{\sI\sW_p,\G_{\mathrm{m},\k}}((\Gr_G)^{\mu_p},\F).
\]
One can also consider the Smith category $\mathrm{Sm}_{\sI\sW_p}((\Gr_G)^{\mu_p},\F)$, defined as a Verdier quotient of this category as in~\eqref{eqn:Smith-cat}, and the corresponding functor $i^{!*}$. We will denote by
\[
\mathsf{Q} \colon \Db_{\sI\sW_p,\G_{\mathrm{m},\k}}((\Gr_G)^{\mu_p},\F) \rightarrow \mathrm{Sm}_{\sI\sW_p}((\Gr_G)^{\mu_p},\F)
\]
the quotient functor.

For any $\lambda\in -\overline{\mathbf{a}_0}\cap \X^\vee$, we have the subcategory 
\[
\Db_{\sI\sW_p,\G_{\mathrm{m},\k}}(\Gr_G^{(\lambda)},\F)\subset \Db_{\sI\sW_p,\G_{\mathrm{m},\k}}((\Gr_G)^{\mu_p},\F)
\]
of complexes supported on $\Gr_G^{(\lambda)}$. Now $\Gr_G^{(\lambda)}$ identifies with the connected component of the base point in a partial affine flag variety, see Theorem~\ref{thm:fixed-pts-Gr}. In $\Db_{\sI\sW_p,\G_{\mathrm{m},\k}}(\Gr_G^{(\lambda)},\F)$ we have a notion of parity complexes, defined by the obvious analogue of Definition~\ref{def:parity}, and a classification theorem similar to Theorem~\ref{thm:classification-parity-1}.
In particular, for any $w\in W^{(\lambda-\xi)}_\af$ we have a ``normalized'' indecomposable parity object $\scE^{\sI\sW_p}_{\lambda,w} \in \Db_{\sI\sW_p,\G_{\mathrm{m},\k}}(\Gr_G^{(\lambda)},\F)$ corresponding to the orbit labelled by $w$. 

Consider now the composition
\begin{multline*}
\Phi\colon \Perv_{\sI\sW}(\Gr_G,\F)\simeq \Perv_{\sI\sW,\G_{\mathrm{m},\k}}(\Gr_G,\F) \\
\hookrightarrow 
\Db_{\sI\sW_p,\G_{\mathrm{m},\k}}(\Gr_G,\F) \xrightarrow{i^{!*}} \mathrm{Sm}_{\sI\sW_p}((\Gr_G)^{\mu_p},\F).
\end{multline*}
The following result is the main geometric result from~\cite{RW22}.

\begin{thm}
\phantomsection
\label{thm:RW}
\begin{enumerate}
\item 
\label{it:RW-Phi-ff}
The restriction of $\Phi$ to tilting perverse sheaves is fully-faithful. 
\item 
\label{it:RW-Phi-image}
For any $\lambda\in -\overline{\mathbf{a}_0}\cap \X^\vee$ and $w\in W^{(\lambda-\xi)}_\af$, we have an isomorphism 
\[
\mathsf{Q}(\scE^{\sI\sW_p}_{\lambda,w})\simeq \Phi(\scT^{\sI\sW}_{w\square_p\lambda}).
\]
\end{enumerate}
\end{thm}

The main ingredient of the proof of Theorem~\ref{thm:RW}\eqref{it:RW-Phi-ff} is the observation that
each $\scT^{\sI\sW}_\mu$ is a parity complex (in the sense that, once again, it satisfies the obvious analogue of the conditions in Definition~\ref{def:parity}), because the dimension of ``relevant'' $I^+_{\un}$-orbits, i.e.~those which support a nonzero Iwahori--Whittaker local system, is of constant parity on each connected component of $\Gr_G$. For the proof of Theorem~\ref{thm:RW}\eqref{it:RW-Phi-image}, ones observes that
we have a theory of parity complexes in $\mathrm{Sm}_{\sI\sW_p}((\Gr_G)^{\mu_p},\F)$ (simply by analysis of the Smith category of a point), with a classification theorem similar to Theorem~\ref{thm:classification-parity-1}. We also use the fact that
the $\G_{\mathrm{m},\k}$-action on $(\Gr_G)^{\mu_p}$ factors through the morphism $t\mapsto t^p$, hence ``looks like the trivial action for coefficient of characteristic $p$."

\subsubsection{Applications in representation theory}

The first application of Theorem~\ref{thm:RW} given in~\cite{RW22} is to a new proof of the linkage principle (Theorem~\ref{thm:linkage}). In fact, 
for $\lambda,\mu\in \X^\vee_{++}$, we have 
\[
\Hom_{\Perv_{\sI\sW}(\Gr_G,\F)}(\sT^{\sI\sW}_\lambda, \sT^{\sI\sW}_\mu)\neq 0 \quad 
\Rightarrow \quad W_\af \square_p \lambda =W_\af \square_p \mu,
\]
because if $W_\af \square_p \lambda \neq W_\af \square_p \mu$ the objects $\Phi(\sT^{\sI\sW}_\lambda)$ and $\Phi(\sT^{\sI\sW}_\mu)$ are supported on distinct connected components of $(\Gr_G)^{\mu_p}$, hence there cannot exist a nonzero morphism between them. By fully faithfulness of $\Phi$ (see Theorem~\ref{thm:RW}\eqref{it:RW-Phi-ff}), this implies that
\[
\Hom_{\Perv_{\sI\sW}(\Gr_G,\F)}(\sT^{\sI\sW}_\lambda, \sT^{\sI\sW}_\mu)=0,
\]
as stated above.
Translated in $\Rep(G^\vee_\F)$, it is not difficult to check that this statement is equivalent to Theorem~\ref{thm:linkage}.

The second application of Theorem~\ref{thm:RW} is a multiplicity formula for indecomposable tilting modules. 
We fix $\lambda\in \overline{C}\cap \X^\vee$; then the elemnt $\lambda+\xi$ is in $-\overline{\mathbf{a}_0}\cap \X^\vee$.

\begin{thm}
\label{thm:tilting-formula}
For any $w,y\in W_\af^{(\lambda)}$, we have 
\[
(\til(w\bullet \lambda):\coweyl(y\bullet \lambda))={}^p \hspace{-1pt} n_{y,w}(1).
\]
\end{thm}

The proof of this theorem proceeds as follows.
It is easily seen that, to check this formula, it suffices to prove that 
\[
\dim \Hom(\til(w\bullet \lambda), \til(y\bullet \lambda))=
\sum_{z\in W_\af^{(\lambda)}} {}^p \hspace{-1pt} n_{z,y}(1) \cdot {}^p \hspace{-1pt} n_{z,w}(1)
\]
for all $w,y\in W_\af^{(\lambda)}$. 
This identity follows from (a slightly more precise version of) Theorem~\ref{thm:RW}\eqref{it:RW-Phi-image} together with the following facts:
\begin{itemize}
\item for any facet $\mathbf{f} \subset \overline{\mathbf{a}_0}$, the pullback functor under the projection $\Fl_{G,p} = \GGr_{P_{\mathbf{a}_0}}\rightarrow \GGr_{P_{f}}$ sends indecomposable parity complexes to indecomposable parity complexes; 
\item the stalks of the indecomposable Iwahori--Whittaker parity complexes on $\Fl_{G,p}$ are described by the polynomials $({}^p \hspace{-1pt} n_{y,w} : y,w \in W_\af)$.
\end{itemize}

\begin{rmk} 
\begin{enumerate}
\item In the case $\lambda=0$, and assuming $p>h$, the formula in Theorem~\ref{thm:tilting-formula} was conjectured in~\cite{RW18}, drawing inspiration from earlier conjecture of Andersen (involving ordinary Kazhdan--Lusztig polynomials). 
In this special case, this formula was proved 
in~\cite{AMRW19}. 
\item The arguments sketched above use a relation between parity complexes on $\GGr_{P_{\mathbf{a}_0}}$ and $\GGr_{P_{\mathbf{f}}}$.  On the other hand, in representation theory, it may happen that $C\cap \X^\vee=\varnothing$ if $p$ is small, so that there is no ``block'' corresponding to the alcove $\mathbf{a}_0$.
\end{enumerate}
\end{rmk}

\subsection{Discussions of the tilting character formula}

\subsubsection{Relation with Lusztig's character formula} 
\label{sss:Lusztigs-formula}

We assume $p\geq h$. The following is the conjecture by Lusztig alluded to in~\S\ref{sss:characters}.

\begin{conjecture}
For any $w\in W^{(0)}_\af$ such that $\langle w\bullet 0+\rho, \alpha\rangle \leq p(p-h+2)$ for all $\alpha\in R_+$, we have 
\[
\mathrm{ch}(\simp(w\bullet 0))=\sum_{y\in W^{(0)}_\af} (-1)^{\ell(w)+\ell(y)} h_{w_0y,w_0w}(1)\cdot \mathrm{ch}(\coweyl(y\bullet 0)).
\]
\end{conjecture} 

As explained in~\S\ref{sss:characters},
this conjecture is known to be true for large $p$.
Using the work of Andersen evoked in~\S\ref{sss:characters}, one can deduce from Theorem~\ref{thm:tilting-formula} another proof of this fact. 

\subsubsection{Relation with Finkelberg--Mirkovi\'{c} conjecture} 

Fix a reductive algebraic group $\mathbf{G}$ over $\F$ whose Frobenius twist $\mathbf{G}^{(1)}$ is $G^\vee_\F$, and assume that $p>h$. 
Let $W_\mathrm{ext}=W\ltimes \X^\vee$ be the extended affine Weyl group. 
We have a decomposition similar to~\eqref{eqn:decomposition-Rep-blocks} for the group $\mathbf{G}$, and we will denote by $\Rep_{\langle 0 \rangle}(\mathbf{G})$ the direct sum of the blocks associated with the elements of the form $w \bullet 0$ where $w$ runs over the elements of $W_\mathrm{ext}$ of length $0$.
We can consider the functor given by pullback via the Frobenius map:
\[
\mathrm{Fr}^*\colon \Rep(G^\vee_\F)\rightarrow \Rep(\mathbf{G}),
\]
and the functor of tensoring with a module over the form $\mathrm{Fr}^*(M)$ with $M \in \Rep(G^\vee_\F)$ stabilizes $\Rep_{\langle 0 \rangle}(\mathbf{G})$.

Consider also the ``opposite affine Grassmannian'' $\Gr_G^{\mathrm{op}} = (\Loop^+G\backslash \Loop G)_{\mathrm{\acute{e}t}}$, the action of $I_{\un}$ on this ind-scheme induced by right multiplication on $\Loop G$, and the associated category of equivariant perverse sheaves. We have a natural bifunctor
\[
(-) \star (-) : \Perv_{\Loop^+G}(\Gr_G, \F) \times \Perv_{I_{\un}} ( \Gr_G^{\mathrm{op}}, \F ) \to \Perv_{I_{\un}}( \Gr_G^{\mathrm{op}}, \F )
\]
which defines a left action of the monoidal category $(\Perv_{\Loop^+G}(\Gr_G, \F), \star)$ on the category $\Perv_{I_{\un}} ( \Gr_G^{\mathrm{op}}, \F)$.

The following conjecture is due to Finkelberg--Mirkovi\'{c}.

\begin{conjecture}
\label{conj:FM}
There exists an equivalence of categories 
\[
\mathrm{FM}\colon \Perv_{I_{\un}} ( \Gr_G^{\mathrm{op}}, \F ) \simto \Rep_{\langle 0 \rangle}(\mathbf{G})
\]
with a bifunctorial isomorphism 
\[
\mathrm{FM}(\sG \star \sF)\simeq \mathrm{FM}(\sF)\otimes_\F \mathrm{Fr}^*(\mathrm{Sat}(\sG))
\]
for any $\scG \in \Perv_{\Loop^+G}(\Gr_G, \F)$ and $\scF\in \Perv_{I_{\un}} ( \Gr_G^{\mathrm{op}}, \F )$. 
\end{conjecture}

\begin{rmk}
\begin{enumerate}
\item
For a discussion of this conjecture and its use in representation theory, see e.g.~\cite{ciappara-williamson}.
\item
For the same reason as for $\Perv_{\sI\sW}(\Gr_G,\F)$ (see~\S\ref{sss:IW-model}), $\Perv_{I_{\un}} ( \Gr_G^{\mathrm{op}}, \F )$ admits a canonical structure of highest weight category. The equivalence in Conjecture~\ref{conj:FM} is expected to intertwine this highest weight structure with the natural one on $\Rep_{\langle 0 \rangle}(\mathbf{G})$ (inherited from the ambient highest weight category $\Rep(\mathbf{G})$).
\item
A proof of this conjecture (under slightly stronger assumptions on $p$) has been recently obtained in joint work with Bezrukavnikov, see~\cite{bez-ri}.
\item
One can also formulate ``singular variants'' of the Finkelberg--Mirkovi\'{c} conjecture that describe geometrically the singular blocks of $\Rep(\mathbf{G})$; see~\cite[Chap.~6, \S 3.2.2]{book} for details.
\end{enumerate}
\end{rmk}

The Finkelberg--Mirkovi\'{c} conjecture (Conjecture~\ref{conj:FM}) and the tilting character formula (Theorem~\ref{thm:tilting-formula}) are philosophically related by Koszul duality:
\[
\begin{tikzcd}
\Db_{I_{\un}} ( \Gr_G^{\mathrm{op}}, \F ) 
\arrow[rr,leftrightarrow,dotted,"\text{Koszul duality}"] \arrow[dr,"\text{Finkelberg--Mirkovi\'{c} conjecture}"'] & 
& \Db_{\sI\sW}(\Fl_G,\F) \arrow[dl,"\text{tilting character formula}"] \\ 
 & \Rep_{\langle 0 \rangle}(\mathbf{G}) &  
\end{tikzcd}
\]
Namely, Theorem~\ref{thm:tilting-formula} expresses the combinatorics of tilting modules in $\Rep_{\langle 0 \rangle}(\mathbf{G})$ in terms of the polynomials ${}^p \hspace{-1pt} n_{y,w}$, which compute dimensions of stalks of parity complexes in $\Db_{\sI\sW}(\Fl_G,\F)$ (see~\cite[\S 11.7]{RW18}). On the other hand the equivalence in the Finkelberg--Mirkovi\'{c} conjecture is expected to be an equivalence of highest weight categories, hence to relate tilting objects in $\Rep_{\langle 0 \rangle}(\mathbf{G})$ to tilting objects in $\Perv_{I_{\un}} ( \Gr_G^{\mathrm{op}}, \F )$. In~\cite{AMRW19}, following ideas of Bezrukavnikov--Yun for characteristic-$0$ coefficients, a ``Koszul duality'' is constructed, which relates parity complexes in $\Db_{\sI\sW}(\Fl_G,\F)$ to tilting objects in a ``graded version'' of $\Perv_{I_{\un}} ( \Gr_G^{\mathrm{op}}, \F )$. 

\bibliographystyle{plain} 
\bibliography{MyBibtexNote}

\end{document}